\documentclass[oneside, 12pt]{amsart}
\usepackage{amscd, amssymb, amsmath, mathrsfs}
\usepackage[english]{babel}
\usepackage{booktabs}
\usepackage{tikz-cd}
\usepackage{url}
\usepackage[pdftex, colorlinks=true,  citecolor=blue, linkcolor=blue, linktocpage=true]{hyperref}

\newcommand\mylabel[1]{\label{#1}\marginpar{\vspace{-1ex}\medskip\medskip\footnotesize \tt #1}}
\renewcommand\mylabel[1]{\label{#1}}
\newcommand{\mydate}{
\number\day\space
\ifcase\month \or January\or February\or March\or April\or May\or June\or July\or August\or September\or October\or November\or December\fi 
\space\number\year}

\DeclareUrlCommand\arXiv{\urlstyle{same}}

\newtheorem{theorem}{Theorem}[section]
\newtheorem*{maintheorem}{Theorem}
\newtheorem{lemma}[theorem]{Lemma}
\newtheorem{proposition}[theorem]{Proposition}
\newtheorem{corollary}[theorem]{Corollary}

\theoremstyle{definition}

\newtheorem*{acknowledgement}{Acknowledgement}

\theoremstyle{remark}

\newcommand{\ZZ}{\mathbb{Z}}

\newcommand{\PP}{\mathbb{P}}

\newcommand{\GG}{\mathbb{G}}

\newcommand{\idealp}{\mathfrak{p}}

\newcommand{\ideala}{\mathfrak{a}}
\newcommand{\idealb}{\mathfrak{b}}

\newcommand{\shI}{\mathscr{I}}
\newcommand{\shJ}{\mathscr{J}}

\newcommand{\shN}{\mathscr{N}}
\newcommand{\shL}{\mathscr{L}}

\newcommand{\alg}{\text{\rm alg}}

\newcommand{\Aut}{\operatorname{Aut}}

\newcommand{\Bl}{\operatorname{Bl}}

\newcommand{\Der}{\operatorname{Der}}

\newcommand{\edim}{\operatorname{edim}}

\newcommand{\Ext}{\operatorname{Ext}}

\newcommand{\Fitt}{\operatorname{Fitt}}

\newcommand  {\Grass}{\operatorname{Grass}}

\newcommand{\Image}{\operatorname{Im}}

\newcommand{\Kernel}{\operatorname{Ker}}

\newcommand{\length}{\operatorname{length}}
\newcommand{\Lie}{\operatorname{Lie}}

\newcommand{\invlim}{\varprojlim}

\newcommand{\lra}{\longrightarrow}

\newcommand{\Mat}{\operatorname{Mat}}

\newcommand{\maxid}{\mathfrak{m}}

\newcommand{\primid}{\mathfrak{p}}
\renewcommand{\O}{\mathscr{O}}

\newcommand{\pdeg}{\operatorname{pdeg}}

\newcommand{\perf}{{\text{\rm perf}}}

\newcommand{\pr}{\operatorname{pr}}
\newcommand{\Proj}{\operatorname{Proj}}

\newcommand{\quadand}{\quad\text{and}\quad}
\newcommand{\quador}{\quad\text{or}\quad}

\newcommand{\ra}{\rightarrow}

\newcommand{\rank}{\operatorname{rank}}
\newcommand{\red}{{\operatorname{red}}}
\newcommand{\Reg}{\operatorname{Reg}}

\newcommand{\Sing}{\operatorname{Sing}}

\newcommand{\Spec}{\operatorname{Spec}}

\newcommand{\Sym}{\operatorname{Sym}}

\newcommand{\trdeg}{\operatorname{trdeg}}

\newcommand{\uHom}{\underline{\operatorname{Hom}}}

\newcommand{\lieg}{\mathfrak{g}}

\begin{document}

\title[Regular curves of genus one]
      {The structure of regular genus-one curves over imperfect fields}

\author[Stefan Schr\"oer]{Stefan Schr\"oer}
\address{Mathematisches Institut, Heinrich-Heine-Universit\"at, 40204 D\"usseldorf, Germany}
\curraddr{}
\email{schroeer@math.uni-duesseldorf.de}


\dedicatory{6 November 2022}

\begin{abstract}
Working over imperfect fields, we give a comprehensive classification of genus-one curves 
that are regular but not geometrically regular, extending the known case of 
geometrically reduced curves. The description is given intrinsically, in terms of
twisted forms of standard models with respect to infinitesimal group scheme actions, and not via extrinsic equations.
The main new idea is to analyze and exploit
moduli for  fields of representatives in Cohen's Structure Theorem.
The results serve for    the understanding
of genus-one fibrations over higher-dimensional bases.
\end{abstract}

\maketitle
\tableofcontents

\section*{Introduction}
\mylabel{Introduction}

Algebraic curves fall into three classes of rather unequal size, each determined by the ampleness properties of the dualizing sheaf: 
The projective line,
elliptic curves, and curves of higher genus. This trichotomy is expected to hold true in   dimensions $d\geq 2$ in the following sense:
Up to passing to  suitable birational models, a proper normal scheme 
$Z$ is either covered by projective lines, or admits a fibration $f:Z\ra B$ over a lower-dimensional base where the generic fiber  satisfies
$c_1=0$, or is of general type. This is indeed true for algebraic surfaces, in light of the
Enriques classification.

To understand the case of a fibration $f:Z\ra B$ whose generic fiber $Y=f^{-1}(\eta)$
satisfies the condition $c_1=0$, one inevitably  is lead to study proper normal schemes $Y$ of dimension $n\geq 1$ over function fields  $F$
such that the structure sheaf satisfies $h^0(\O_Y)=1$ and the dualizing sheaf  $\omega_Y$ is numerically trivial.
The goal of this paper is to carry out an in-depth analysis for the simplest relevant case, when
$Y$ is a proper normal curve with numerical invariants $h^0(\O_Y)=h^1(\O_Y)=1$ and dualizing sheaf $\omega_Y=\O_Y$.
We also say that $Y$ is a \emph{genus-one curve}.

The situation is   completely understood in characteristic zero, when such an $Y$ is a principal homogeneous space
over an elliptic curve. I like to say that $Y$ is \emph{para-elliptic}, compare \cite{Laurent; Schroeer 2021} and \cite{Schroeer 2022b}.
One also understands the case when $Y$ is geometrically singular and geometrically reduced: 
Then the curve is  a twisted form of the rational cuspidal curve $\Spec k[u^2,u^3]\cup \Spec k[u^{-1}]$,
and this happens precisely in  characteristic is $p\leq 3$. 
The occurring equations where analyzed  by Queen (\cite{Queen 1971}, \cite{Queen 1972}),  and one also says that $Y$ is \emph{quasi-elliptic}.
Note that the occurrence of such curves is  arguably the key phenomenon in the extension of the Enriques classification
to positive characteristics by Bombieri and Mumford (\cite{Bombieri; Mumford 1976}, \cite{Bombieri; Mumford 1977}).
In arithmetic settings, the understanding  of quasi-elliptic fibrations at the primes two or three is often a key step
to unravel the situation over fields of characteristic zero (for example in \cite{Schroeer 2023}).

Throughout this paper, we investigate the general  case, \emph{when the field $F$ is  imperfect  of characteristic $p>0$,
and $Y$ is  regular but not geometrically regular}. Note that this includes the geometrically non-reduced curves,
and our results mainly pertain to them.
I already touched upon the topic in \cite{Schroeer 2010}, Section 6, but with  very limited results.

Here we use a completely different and novel approach, which relies on a connection between commutative algebra
and algebraic geometry that to my best knowledge was not explored so far:
The \emph{fields of representatives} occurring in Cohen's Structure Theorems for complete local noetherian rings
are usually non-unique, and thus have \emph{moduli  in the sense of algebraic geometry}.
The main idea of this paper  is to study and exploit such moduli in a systematic fashion. I expect that this  
should be useful   in many other contexts as well.
 
To state our main  result, we have to introduce certain curves  
$C=C^{(i)}_{r,F,\Lambda}$ called \emph{standard models}.
They are birational to the projective line $\PP^1_R$ over the local Artin ring $R=F[W_1,\ldots,W_r]/(W_1^p,\ldots,W_r^p)$,
in which the field $F$ and the integer  $r\geq 0$ enter.
The birational morphisms $\PP^1_R\ra C$ is specified in terms of certain $F$-subalgebras $\Lambda$ lying inside the ring of dual numbers $R[\epsilon]$ or the product ring $R\times R$,
subject to   conditions that among other things ensure $h^0(\O_C)=h^1(\O_C)=1$. The  upper index pertains to this subring,
and is  a symbol that could take  the values $0\leq i\leq r$ or $i=(1,1)$.
Note that  one views the ring of dual numbers or the product ring as coordinate rings of some effective Cartier divisor on $\PP^1_R$,
and that the subring $\Lambda$ has moduli.
For details of the construction we refer to  Section \ref{Formulation}.  We now can state our first main result:

\begin{maintheorem}
(see Thm.\ \ref{properties genus-one})
Let $Y$ be a genus-one curve that is regular but not geometrically regular. Then the following holds:
\begin{enumerate}
\item The field $F$ has characteristic  at most three,  and $p$-degree at least $ r+1$.
\item In characteristic two the second Frobenius base-change   $Y^{(p^2)}$ is a twisted form of some standard model $C^{(i)}_{r,F,\Lambda}$
with $0\leq i\leq 2$ or $i=(1,1)$.
\item In characteristic three, the first Frobenius base-change $Y^{(p)}$ is isomorphic
to some standard model $C^{(i)}_{r,F,\Lambda}$ with  $0\leq i\leq 1$.
\end{enumerate}
\end{maintheorem}

In particular, the base-change $Y\otimes F^\alg$ is one of the specified standard models over the algebraic closure $F^\alg$.
Concerning   existence, we have the following:

\begin{maintheorem}
(see Thm.\ \ref{existence genus-one})
Let  $r\geq 0$ be some integer. Suppose   characteristic and   $p$-degree of the ground field $F$, and the symbol $i$
are as in one of the rows of the following table:
$$
\begin{array}{llll}
\toprule
\operatorname{char}(F)	& \pdeg(F)	& i	\\
\midrule
3		& \geq r+1	& 0,1 	\\
2		& \geq r+1	& 0,1,2	\\
2		& \geq r+2	& (1,1)	\\
\bottomrule
\end{array}
$$
Then for some $\Lambda$ the standard model $C^{(i)}_{r,\Lambda,F}$ admits a   twisted form $Y$ that is regular.
\end{maintheorem}

I find these results surprising on at least two counts: First, due to   moduli for the subalgebras $\Lambda$ 
inside $R[\epsilon]$ or $R\times R$,   the genus-one curves $Y$ are \emph{usually not 
twisted forms    of some standard model without making a ground field extension}. This is in marked contrast two the case of
para-elliptic curves, which are twisted forms of their Jacobian, and also for    quasi-elliptic curves,
all of which are twisted forms of the rational cuspidal curve.
Second, the existence of genus-one curves that are regular but   geometrically non-reduced remains  \emph{confined to characteristic two and three}:
This is as with   quasi-elliptic curves, a state of affairs that I did not expect.

Let me point out that our approach emphasizes intrinsic aspects of the curves,
and also sheds further light on their geometry. I hope that
future research will   lead to   equations as well. However, I expect that 
such equations to be rather complicated and of extrinsic nature, and thus perhaps of limited practical value.  

The idea for the proof for the above theorems is as follows: Using a result from \cite{Schroeer 2010},
the problem is reduced to the case that for some simple height-one extension  $F'=F(\omega^{1/p})$,
the base-change $Y'=Y\otimes_FF'$ becomes singular, yet stays integral. We then argue that the normalization
$X'$ must be a genus-zero curve, over some further height-one extension $F'\subset E'$, a phenomenon called \emph{constant field extension}
in classical parlance.
Using the comparatively simple structure of regular genus-zero curves, combined with the arithmetic of the height-one extensions,
we can unravel the geometry of the \emph{conductor square} that accompanies the normalization map.
It is  precisely this interplay that yields the connection to the standard forms $C_{r,F,\Lambda}^{(i)}$.
We also get an  explicit description of the sheaf $\Omega^1_{Y'/F'}$, which has to be locally free modulo 
its torsion part.  Using the techniques already developed in  \cite{Schroeer 2008} and \cite{Brion; Schroeer 2022},
this allows to discards all primes except $p\leq 3$.

Any classification of proper smooth schemes $X$ with $c_1=0$ of dimension $n\geq 3$ over an algebraically closed ground field $k$
of positive characteristics is likely to involve the above classification, in parallel to  the Enriques classification
by Bombieri and Mumford. More general, it should be important for the
understanding fibrations $f:X\ra B$, of relative dimension one or higher, whose generic fiber has $c_1=0$, as appearing in
the Minimal Model Program and  Mori Theory.
One could also hope to generalize the  Ogg--Shafarevich Theory or the description of multiple fibers
(compare \cite{Dolgachev; Gross 1994}, \cite{Gross 1997}, \cite{Szydlo 2004}).

\medskip
The paper is organized as follows:
In Section \ref{Generalities} we collect some foundational facts on denormalization and twisting that are used throughout.
The main results of the paper are formulated in Section \ref{Formulation}.
Section \ref{Intersection algebras} establishes  some relevant facts on intersections of subrings in Artin rings.
We then start to construct regular genus-one curves that are not geometrically regular:
In Section \ref{Constructions subfields} we use pairs of subfields, a particularly simple way with direct geometric appeal,
which however works only in characteristic two.
Section \ref{Constructions fields representatives} contains our construction relying on non-standard field of representatives. This also involves 
the choice of certain additive vector fields. An analysis of the restricted Lie algebra of all global vector fields
is given in Section \ref{Lie algebra}.
Section \ref{Constructions nilpotents} contains the most challenging  construction of genus-one curves, where nilpotent
elements play a crucial role in the denormalization. This relies on some technical observations from commutative algebra 
verified  in Section \ref{Verification}.
After Section \ref{Twisted ribbons}, where we collect some results on  ribbons and their twisted forms
 in connection to genus-zero curves, we examine genus-one curves with singularities
in Section \ref{Genus-one curves}. The paper culminates in Section \ref{Proof}, which comprises the proofs for our main results.

\begin{acknowledgement}
The research was   conducted       in the framework of the   research training group
\emph{GRK 2240: Algebro-Geometric Methods in Algebra, Arithmetic and Topology}.
\end{acknowledgement}

\section{Generalities and recollections}
\mylabel{Generalities}

Let $F$ be a ground field, for the moment  of arbitrary characteristic $p\geq 0$.
Throughout, the term \emph{curve} refers to a proper scheme $X$ that is equi-dimensional,
of dimension one. It comes with the numerical invariants $h^i(\O_X)=\dim_F H^i(X,\O_X)$,
for  $0\leq i\leq 1$.  
We say that $X$ is a \emph{curve of genus} $g$ if these invariants take the values $h^0(\O_X)=1$ and $h^1(\O_X)=g$.
Clearly, the condition is stable under ground field extensions.
Moreover, such curves are geometrically connected and without embedded components, hence Cohen--Macaulay.  
Note that we do not impose any other conditions whatsoever on the singularities, and in particular allow nilpotent elements
in the structure sheaf.

For $g=0$ we   say that $X$ is a \emph{genus-zero curve}. Examples are, of course, the projective line,
or quadric curves in $\PP^2$, or   infinitesimal extensions $X=\PP^1\oplus\shI$, where $\shI$ 
is a locally free sheaf on $\PP^1$ with splitting type $(-1,\ldots,-1)$.
For $g=1$ we say that  $X$ is a  \emph{genus-one curve}. Here   examples are elliptic curves,
or cubic curves in $\PP^2$, or   infinitesimal extensions $X=E\oplus \shI$, where $E$ is an elliptic curve and  $\shI$ is 
a sum of invertible sheaves that are non-trivial of degree zero, or $X=\PP^1\oplus \shJ$, where $\shJ$ has splitting type
$(-2,-1,\ldots,-1)$.

Let $X$ be an arbitrary curve, and assume for simplicity that $X$ is irreducible and has no embedded components.
Let $f:X\ra Y$ be  a finite birational morphism between such curves.
The \emph{branch scheme} $B\subset Y$ is the closed subscheme corresponding to the \emph{sheaf of conductor ideals} $\shI\subset\O_Y$,
which is defined as the annihilator   of $f_*(\O_X)/\O_Y$.
The preimage $A=f^{-1}(B)$ is called the \emph{ramification scheme}. It is the closed subscheme of $Y$
corresponding to  sheaf of ideals $\shI\subset f_*(\O_X)$. Indeed, by the very definition of $\shI\subset\O_Y$,
this is also a sheaf of ideals inside $f_*(\O_X)$. 
Moreover, the resulting cartesian diagram
\begin{equation}
\label{conductor square}
\begin{CD}
A	@>>>	X\\
@VVV		@VVfV\\
B	@>>>	Y
\end{CD}
\end{equation}
is also cocartesian. Then $0\ra\O_Y\ra f_*(\O_X)\oplus\O_B\ra f_*(\O_A)\ra 0$ is exact,
where the first arrow is the diagonal map, and the second arrow is the difference map.
In turn, we have a long exact sequence
\begin{equation}
\label{conductor sequence}
0\ra \Gamma(\O_Y)\ra\Gamma(\O_X)\oplus\Gamma(\O_B)\ra\Gamma(\O_A)\ra H^1(Y,\O_Y)\ra H^1(X,\O_X)\ra 0.
\end{equation}
Note that both maps $\O_Y\ra f_*(\O_X)$ and $\O_B\ra f_*(\O_A)$ are injective. Let us record the following fundamental facts:

\begin{proposition}
\mylabel{fundamental facts}
In the above situation, the following holds:
\begin{enumerate}
\item
The numerical invariants of the schemes are related by the formula
$$
h^0(\O_Y) +h^0(\O_A)  + h^1(\O_X) = h^0(\O_X)+h^0(\O_B) +h^1(\O_Y).
$$
\item 
Suppose the restriction map $\Gamma(\O_X)\ra\Gamma(\O_A)$ is injective. Then the ring $\Gamma(\O_Y)$
equals the intersection    $\Gamma(\O_X)\cap \Gamma(\O_B)$ inside $\Gamma(\O_A)$.
\item 
The finite morphism $f:X\ra Y$ is a   universal homeomorphism if and only if this holds for $f:A\ra B$.
\item 
If   the one-dimensional scheme $Y$ is Gorenstein then  for each point $b\in B$
the formula $\sum_{a\in f^{-1}(b)}h^0(\O_{A,a}) = 2h^0(\O_{B,b})$ holds. The converse   is true
provided that $X$ is Gorenstein and $A\subset X$ is Cartier.
\end{enumerate}
\end{proposition}

\proof
Statement (i) follows from the additivity of Euler characteristics and the exact sequence \eqref{conductor sequence}.
The latter also gives (ii). Assertion (iii) is a consequence of \cite{EGA IVd}, Corollary 18.12.11,
and (iv) follows from \cite{Fanelli; Schroeer 2020}, Proposition A.2 and A.3, compare also the references discussed there.
\qed

\medskip
We will use the following consequence frequently, where here and throughout the symbol $\epsilon$ denotes an indeterminate
subject to the condition $\epsilon^2=0$:

\begin{corollary}
\mylabel{subring}
Suppose $Y$ is  Gorenstein and $f:X\ra Y$ is a universal homeomorphism.
Let $a\in A$ be a closed point,  with image $b\in B$ and residue field $E$. Assume that $\O_{A,a}=E[\epsilon]$ is a ring of dual numbers.
Then the subring $\O_{B,b}$ takes the form
$\Lambda=L+H\epsilon$, where $L\subset E[\epsilon]$ is a subfield,
and $H\epsilon\subset E\epsilon$ is an $L$-linear subspace of codimension one.
\end{corollary}

\proof
Write   $L=\kappa(b)$ for the residue field, such that we have extensions $F\subset L\subset E$.
By Cohen's Structure Theorem (\cite{AC 8-9}, Chapter IX, \S4, No.\ 3, Theorem 1), the residue class map $\O_{B,b}\ra L$ admits a section, giving
an inclusion of  $L$ into $\O_{B,b}$ and $\O_{A,a}$. Write $\maxid_b$ and $\maxid_a$ for the
maximal ideals in these local Artin rings.
The inclusion $\maxid_b\subset\maxid_a=E\epsilon$   shows
that $\maxid_b^2=0$, hence the $\O_{B,b}$-module structures comes from a structure of    $L$-vector spaces.  
Now write  $\maxid_b=H\epsilon$. Clearly
$$
h^0(\O_{A,a}) = 2[E:L]\cdot [L:F]\quadand h^0(\O_{B,b}) = \left(1+\dim_L(H\epsilon)\right) \cdot [L:F].
$$
According to the proposition, we have
$h^0(\O_{A,a})=2h^0(\O_{B,b})$. Combining with the above equations, we arrive at
$[E:L]=1+\dim_L(H\epsilon)$, and thus $H\epsilon\subset E\epsilon$ must be an $L$-hyperplane.
\qed

\medskip
Note that the sum $\Lambda=L+H\epsilon$ is direct because $L\cap H\epsilon=0$. However, this sum is usually
not compatible with the canonical decomposition $E[\epsilon]=E\oplus E\epsilon$.
Also note that the local Artin ring $\O_{B,b}$ is a field if and only if $H\epsilon=0$, which in turn means
that the subfield $L\subset\O_{A,a}$ is a field of representatives.
 
Recall that to form    diagram \eqref{conductor square}, we 
started with a finite birational morphism $f:X\ra Y$, which then determines $B$ and $A=f^{-1}(B)$. 
Conversely, if one begins with a curve $X$, for simplicity assumed to be irreducible and without embedded components,
together with a finite closed subscheme $A\subset X$, a finite scheme $B$, and an inclusion $\Gamma(\O_B)\subset\Gamma(\O_A)$,
one obtains a curve $Y$ by forming the cocartesian diagram \eqref{conductor square}.
The push-out $Y$ indeed exists as an algebraic space (\cite{Artin 1970}, Theorem 6.1), which  here is actually a scheme (\cite{Ferrand 2003}, Theorem 7.1). 
Such constructions are also called \emph{pinchings}.
It is not difficult to characterize those pinchings where  $A$ and $B$ are actually defined by the sheaf of conductor ideals:  

\begin{lemma}
\mylabel{pinching with conductor}
In the above situation, the  diagram \eqref{conductor square} is a conductor square if and only if the schematic support of
$ f_*(\O_A)/\O_B$ coincides with  $B$.
\end{lemma}

\proof
Let $B'\subset B$ be the schematic support of the coherent sheaf $ f_*(\O_A)/\O_B$.
The sequence $0\ra\O_Y\ra f_*(\O_X)\oplus\O_B\ra f_*(\O_A)\ra 0$ is exact, to the left
because $Y$ has no embedded components and $f$ is birational, in the middle 
by  the universal property of cocartesian squares,
and to the right because $A\subset X$ is a closed embedding. Applying \cite{Kashiwara; Schapira 2006}, Lemma 8.3.11 to the resulting
cocartesian square of abelian sheaves on $Y$, we see that the canonical map $f_*(\O_X)/\O_Y\ra f_*(\O_A)/\O_B$ is bijective.
It immediately follows that the condition is necessary. For the converse, suppose that  $B'=B$. It remains to check
that the inclusion $A\subset f^{-1}(B)$ is an equality. Applying loc.\ cit.\ again, we see that the canonical map
$$\Kernel(\O_Y\ra\O_B)\lra\Kernel(f_*(\O_X)\ra f_*(\O_A))$$ is bijective, hence $A=f^{-1}(B)$.
\qed

\medskip
The prime objective of the paper is to understand the structure of   genus-one curves $Y$ that are regular
but not geometrically regular. The following well-known facts will be important:

\begin{proposition}
\mylabel{genus-zero and genus-one}
Each genus-one curve $Y$ that is integral and Gorenstein has dualizing sheaf $\omega_Y=\O_Y$.
Each genus-zero curve $X$    that is integral and Gorenstein is isomorphic to a quadric curve in $\PP^2$.
If it admits an invertible sheaf of degree one, we have $X\simeq \PP^1$.
\end{proposition}

\proof
The invertible sheaf $\shN=\omega_Y$ has degree $d=-2\chi(\O_Y)=0$. It admits a global section $s\neq 0$, because  $h^0(\omega_Y)=h^1(\O_Y)=1$.
The map $s:\O_Y\ra\shN$ is injective because $Y$ is integral, and thus must be bijective since $d=0$.

The invertible sheaf $\shL=\omega_X^{\otimes-1}$ has degree $d= 2\chi(\O_X)=2$.
The usual arguments with Serre Duality and Riemann--Roch show that $\shL$ is globally generated
with $h^0(\shL)=3$,
and we obtain  a morphism $f:X\ra\PP^2$ with $\shL=f^*\O_{\PP^2}(1)$. Let $C\subset \PP^2$ be the image,
which is an integral curve of degree at least two.
The Degree Formula $\deg(\shL)=\deg(f)\cdot \deg(C)$ reveals that $f$ is birational and $C$ is a quadric.
Forming the conductor square \eqref{conductor square} and applying  Proposition \ref{fundamental facts} with $f:X\ra C$, we see that 
$h^0(\O_A)=h^0(\O_B)$, thus $f$ is an isomorphism.

If there is an invertible sheaf $\shN$ of degree one, the usual arguments show that $\shL$ is globally generated
with $h^0(\shN)=2$, and that the resulting morphism $g:X\ra\PP^1$ with $\shN=g^*(\O_{\PP^1}(1)$ is an isomorphisms.
\qed

\medskip
We close  this section with some general observations on \emph{twisting}.
Let $Y$ be a scheme over our ground field $F$. Another scheme $\tilde{Y}$ is called a \emph{twisted form}
if $\tilde{Y}\otimes E\simeq Y\otimes E$ for some field extension $F\subset E$.
In characteristic $p>0$, such twisted forms may arise as follows: 
Write $\Theta_{Y/F}=\uHom(\Omega^1_{Y/F},\O_Y)$ for the tangent sheaf. 
Let $\lieg\subset H^0(X,\Theta_{X/F})$ be a finite-dimensional subspace that is stable
with respect to Lie bracket $[D,D']$ and $p$-map $D^{[p]}$. Then $\lieg$ is a finite-dimensional \emph{restricted Lie algebra},
and thus corresponds to an infinitesimal group scheme $G$ of height one with  $\Lie(G)=\lieg$.
Moreover, the inclusion into $H^0(X,\Theta_{X/F})$  is nothing but a faithful $G$-action on $Y$,
all this by the Demazure--Gabriel correspondence, see \cite{Demazure; Gabriel 1970}, Chapter II, §7, Theorem 3.5.
We also refer to \cite{Roessler; Schroeer 2022}, Section 1 for more background.
If $T=\Spec(E)$ is a $G$-torsor, one may form the twisted form
$$
\tilde{Y}=T\wedge^G Y = G\backslash (T\times Y),
$$
where the $G$-action on the product is given by $\sigma\cdot (t,y)=(\sigma t, \sigma y)$.
Note that the $G$-action is free, so the quotient exists as an algebraic space (for example  \cite{Laurent; Schroeer 2021}, Lemma 1.1),
which here must be a scheme  (\cite{Olsson 2016},  Theorem 6.2.2, compare also \cite{Schroeer; Tziolas 2023}, Lemma 4.1).
For examples of general twists that become non-schematic one may consult \cite{Schroeer 2022a}.
 
We now discuss a criterion that such twisted forms become regular.
Given a quasicoherent sheaf of ideals $\shI_0\subset\O_Y$, 
there is a largest $\lieg$-stable sheaf of ideals $\shI$ contained in $\shI_0$.
It is actually quasicoherent, and the corresponding closed subscheme is the \emph{orbit} $Z=G\cdot Z_0$ for the closed subscheme $Z_0\subset Y$
defined by $\shI_0$.  

\begin{lemma}
\mylabel{twisting removes singularities}
In the above situation, suppose that $Y$ is noetherian. Then the noetherian scheme  $\tilde{Y}$ is regular
if the following three conditions hold:
\begin{enumerate}
\item The twisted form $\tilde{Z}=T\wedge^G Z$ is regular.
\item The sheaf of ideals for $Z\subset Y$ is locally generated by a regular sequence.
\item For the   open set $U=Y\smallsetminus Z$, the base-change $U\otimes E$ is regular.
\end{enumerate}
\end{lemma}

\proof
The projection $U\otimes E\ra U$ is faithfully flat, so $U$ is regular  by \cite{EGA IVb}, Corollary 6.5.2.
Now fix some $b\in \tilde{Y}$ that belongs to the closed subscheme $\tilde{Z}=\tilde{Y}\smallsetminus \tilde{U}$. Then the local ring $\O_{\tilde{Y},b}$
is regular according to  \cite{EGA IVa}, Chapter 0, Corollary 17.1.9.
\qed

\medskip
Note that the above principal already appeared in somewhat different form in in  \cite{Schroeer 2007}, Proposition 2.2.
Also note that the situation becomes particularly simply for the restricted Lie algebra $\lieg=k$, where both bracket and $p$-map are zero.
The inclusion into $H^0(Y,\Theta_{Y/F})$ corresponds to a global vector field $D$ that is non-zero and satisfies $D^{[p]}=0$, and the group scheme
$G=\alpha_p$ is the Frobenius kernel for the additive group $\GG_a$.
Moreover, the  sheaf of ideals $\shI\subset \O_Y$ defining the orbit $Z=G\cdot Z_0$ is the intersection
for the kernels   for the composite  maps $D^i:\shI_0\ra\O_Y/\shI_0$, $1\leq i\leq p-1$.
Likewise is the case $\lieg=\mathfrak{gl}_1(k)$, which stands for the vector space $k$ endowed with
the trivial bracket and   $p$-map  given by $\lambda^{[p]}=\lambda^p$.
Now the inclusion into $H^0(Y,\Theta_{Y/F})$ is given by a non-zero global vector field $D$ satisfying $D^{[p]}=D$.

\section{Formulation of the main results}
\mylabel{Formulation}

In this section we state our main results on the structure of regular genus-one curves.
For this we have to introduce  certain highly singular genus-one curves, the so-called \emph{standard models} $C=C^{(i)}_{r,F,\Lambda}$.
Let $F$ be a ground field of characteristic $p>0$.
Recall that a \emph{$p$-basis}
for a finite  $F$-algebra $R$  is a set of elements $\omega_1,\ldots,\omega_r\in R$ such that $\omega_i^p\in F$, and  that the resulting
map $F[T_1,\ldots,T_r]/(T_1^p-\omega_1^p,\ldots,T_r^p-\omega_r^p)\ra R$ given by $T_i\mapsto \omega_i$ is bijective.
In other words, each member from $R$ can be written in a unique way as a polynomial $P(\omega_1,\ldots,\omega_r)$, with all exponents at most $p-1$
and all coefficients from $F$. In such a setting, we   say that $P(T_1,\ldots,T_r)$ is a  \emph{$p$-truncated polynomial}.
For more on these notions, see \cite{Kiehl; Kunz 1965}, Section 2, compare also \cite{A 4-7}, Chapter V, \S13, No.\ 1.

Now fix some integer $r\geq 0$, and consider  
$R=F[W_1,\ldots,W_r]/(W_1^p,\ldots,W_r^p)$. Clearly, the classes $w_i\in R$ of the indeterminates $W_i$ yield a  $p$-basis.
We now form the relative projective line
$$
\PP^1_R=\Proj R[T_0,T_1] = (\Spec R[u]) \cup (\Spec R[u^{-1}]),
$$
where $u=T_1/T_0$ is the affine coordinate function.

The   standard models $C=C^{(i)}_{r,F,\Lambda}$  are   defined for $0\leq i\leq r$  in terms of the  closed subscheme  $A=V_+(T_1^2)$.
Its coordinate ring $R[u]/(u^2)=R[\epsilon]$ can be seen as a \emph{ring of dual numbers}.
We now consider subrings $\Lambda\subset R[\epsilon]$ of the form $\Lambda=L+ H\epsilon$,
where $L\subset R[\epsilon]$ is a subring admitting a $p$-basis of length $r-i$, and $H\epsilon\subset R\epsilon$
is a  an $L$-submodule such that $R\epsilon/H\epsilon$ free of rank one as $L$-module.
Furthermore, we demand that    $L\cap R\epsilon=0$, and that  the intersection $\Lambda\cap R\subset R[\epsilon]$
coincides with $F$.
Under all these assumptions   we  form the cocartesian square
$$
\begin{CD}
\Spec (R[\epsilon])	@>>>	\PP^1_R\\
@VVV			@VVV\\
\Spec (\Lambda)	@>>>	C,
\end{CD}
$$
which defines the \emph{$i$-th standard model}   $C=C^{(i)}_{r,F,\Lambda}$. 
The following is a first indication   of their relevance for the classification of regular genus-one curves:

\begin{proposition}
\mylabel{standard models}
The  standard models $C^{(i)}_{r,F,\Lambda}$ are geometrically irreducible genus-one curves whose local rings
are Gorenstein and geometrically unibranch.  
\end{proposition}

\proof
It suffices to prove this when the ground field $F$ is algebraically closed.
Set $X=\PP^1_R$ and $Y=C^{(i)}_{r,F,\Lambda}$.
Clearly $X_\red$ is the projective line,  and  the composite map $\PP^1\ra X\ra Y$ is bijective, and the     residue fields extensions
are equalities. Thus $Y$ is irreducible, and the local rings $\O_{Y,y}$ are geometrically unibranch.

Obviously $h^0(\O_X)=p^r$ and $h^1(\O_X)=0$ and   $h^1(\O_A)=2p^r$. So Proposition \ref{fundamental facts} gives the formula
$h^0(\O_Y) + p^r = h^0(\O_B) + h^1(\O_Y)$.
By our assumptions on the subring $\Lambda=L+H\epsilon$, we have   $[L:F]=p^{r-i}$ and $\dim_F(R\epsilon/H\epsilon)= [L:F]=p^{r-i}$.
Using $\dim_F(R\epsilon)=[R:F]=p^r$ we get $h^0(\O_B)=p^{r-i} + (p^r-p^{r-i})=p^r$, and  conclude $h^0(\O_Y)=h^1(\O_Y)$. 
To see that $Y$ is a genus-one curve it thus suffices to verify
$h^0(\O_Y)=1$. 
Clearly, the canonical map $\Gamma(\O_X)\ra\Gamma(\O_A)$ is injective, so Proposition \ref{fundamental facts}
gives $\Gamma(\O_Y)=\Gamma(\O_X)\cap\Gamma(\O_B)$, where the  intersection takes place inside $\Gamma(\O_A)$.
Again by our assumption on $\Lambda$,  the intersection   coincides with $F$, and we get $h^0(\O_Y)=1$. 
This shows that $Y$ is a genus-one curve.

One immediately sees that $h^0(\O_{A,a})=2h^0(\O_{B,b})$ for every $a\in A$, with $b=f(a)$.
By Proposition \ref{fundamental facts}, we see that $Y$ is Gorenstein. 
\qed

\medskip
There is another relevant standard model $C=C^{(1,1)}_{r,F,\Lambda}$. It is defined   in terms of the disconnected
closed subscheme $A=V_+(T_0T_1)$, whose coordinate ring is $R\times R=R[u]/(u)\times R[u^{-1}]/(u^{-1})$.
Consider subrings $\Lambda\subset R\times R$ of the form 
$\Lambda=\Lambda'\times\Lambda''$, where each factor admits a $p$-basis of length $r-1$,
and $\Lambda'\cap \Lambda''\subset R$ coincides with $F$.
The ensuing cocartesian diagram
$$
\begin{CD}
\Spec (R\times R)	@>>>	\PP^1_R\\
@VVV			@VVfV\\
\Spec (\Lambda'\times\Lambda'')	@>>>	C,
\end{CD}
$$
defines the   standard model $C=C^{(1,1)}_{r,F,\Lambda'\times\Lambda''}$. We also designate 
it by $ C^{(i)}_{r,F,\Lambda}$, with upper index $i=(1,1)$ and $\Lambda=\Lambda'\times\Lambda''$.
Note that this curve fails to be Gorenstein for $p\neq 2$, because $h^0(\O_{A,a})=ph^0(\O_{B,b})$, with $a\in A$ and  $b=f(a)$.
However, as above one checks:

\begin{proposition}
\mylabel{standard model 11}
In characteristic two, the standard models $C^{(1,1)}_{r,F,\Lambda'\times\Lambda''}$ are geometrically irreducible genus-one curves whose local rings
are Gorenstein and geometrically unibranch. 
\end{proposition}
 
For each $F$-scheme $Y$ and $t\geq 0$, we write $Y^{(p^t)}=Y\otimes_FF$ for the base-change with respect to the $t$-fold Frobenius map
$\lambda\mapsto \lambda^{p^t}$ on the ground field, and call them the  \emph{ Frobenius base-changes}. 
We are now ready to  formulate the main results of this paper. 

\begin{theorem}
\mylabel{properties genus-one}
Let   $Y$ be a   genus-one curve that is regular but not geometrically regular, and $r=\edim(\O_{Y,\eta}/F)$ be its geometric generic embedding dimension. 
Then the following holds:
\begin{enumerate}
\item The ground field $F$ has characteristic  $p\leq 3$,  and $p$-degree at least $r+1$.
\item In characteristic two, the second Frobenius base-change $Y^{(p^2)}$ is isomorphic to  some standard model $C^{(i)}_{r,F,\Lambda}$ with
$0\leq i\leq 2$ or $i=(1,1)$. 
\item In characteristic three, the first Frobenius base-change $Y^{(p)}$ is isomorphic
to some standard model $C^{(i)}_{r,F,\Lambda}$ with  $0\leq i\leq 1$.
\end{enumerate}
\end{theorem}

Recall that the vector space dimension of $\Omega^1_{F/\ZZ}$ is called    \emph{$p$-degree} of the field $F$, and written
as  $n=\pdeg(F)$.
This   invariant is also determined by the formula $[F:F^p]=p^n$, and serves as  a measure for imperfectness.
Note that $\pdeg(F)=\dim(W)=\trdeg_k(F)$ provided that $F$
is the function field of an integral scheme $W$ of finite type over a perfect field $k$.
The \emph{geometric generic embedding dimension} $r=\edim(\O_{Y,\eta}/F)$
is the embedding dimension of the local Artin ring $\O_{Y,\eta}\otimes_FF^\perf$. This  
was introduced in \cite{Fanelli; Schroeer 2020}, Section 1 as a measure of geometric non-reducedness. The statement  $\pdeg(F)>\edim(\O_{Y,\eta}/F)$ 
in the above theorem is a general fact established in \cite{Schroeer 2010}, Theorem 2.3, included here for the sake of clarity and  completeness.
As to existence, we have:

\begin{theorem}
\mylabel{existence genus-one}
Let  $r\geq 0$ be some integer. Suppose characteristic and   $p$-degree of the ground field $F$, and the symbol $i$
are as in one of the lines of the following table:
$$
\begin{array}{llll}
\toprule
\operatorname{char}(F)	& \pdeg(F)	& i	\\
\midrule
3		& \geq r+1	& 0,1 	\\
2		& \geq r+1	& 0,1,2	\\
2		& \geq r+2	& (1,1)	\\
\bottomrule
\end{array}
$$
Then for some $\Lambda$ the ensuing standard model  $C^{(i)}_{r,\Lambda,F}$ admits  a   twisted form $Y$ that is regular.
\end{theorem}
 
Summing up, the above two results give a complete classification of the regular genus-one curves
that are not geometrically regular. 
The proofs   require  extensive preparation, and will be given in Section \ref{Proof}.

\section{Intersection algebras in    Artin rings}
\mylabel{Intersection  algebras}

In this section we establish some technical results on intersections
of subalgebras, which already played a role in the definition of the standard models
$C=C^{(i)}_{r,\Lambda,F}$.

Let $F$ be a ground field of characteristic $p\geq 0$,
and $R$ be some finite $F$-algebra, endowed with   ideals $\ideala_1,\ldots,\ideala_s$
together with i integers $n_1,\ldots,n_s\geq 0$.
Given a ring  extension $F\subset A$, we use index notation  $R_A=R\otimes A$ for the base-change.
We now consider the \emph{Grassmann schemes} $\Grass^{n_i}_{R/F}$ whose $A$-valued points
are the $A$-submodules $L_i\subset R_A$ whose quotients   are locally free of rank
$[R:F]-n_i$. Note that if  $A=K$ is a field, this boils down to the vector subspaces $L_i\subset R_K$ of dimension $\dim_K(L_i)=n_i$.
We refer to \cite{EGA I}, Section 9.7 for a comprehensive treatment.

\begin{proposition}
\mylabel{first locally closed}
There is a locally closed set $Z\subset\prod_{i=1}^s\Grass^{n_i}_{R/F}$
such that for each field extension $F\subset K$, a  $K$-valued point  $(L_1,\ldots,L_s)$ belongs to $Z$ if and only
if the $L_i\subset R_K$ are   subalgebras, the residue class maps $L_i\ra (R/\ideala_i)_K$ are injective,
and $K=L_1\cap\ldots\cap L_r$. 
\end{proposition}

\proof
First observe that a vector subspace $L\subset R$ is a subalgebra if and only if 
the map $L\otimes L\ra R/L$ induced by multiplication is zero, and also the projection $F\ra R/L$
is zero. These are closed conditions.
Moreover, the residue class map  $L\ra R/\ideala_i$ is injective if an only $\Lambda^n(L)\ra \Lambda^n(R/\ideala_i)$
is non-zero, where $n=[L:F]$. This is an open condition. 
Given  vector subspaces $L_1,\ldots,L_s\subset R$, the condition $F=L_1\cap\ldots\cap L_s$
means  that the diagonal map $\varphi:R\ra R/L_1\times\ldots\times R/L_r$  
vanishes on $F$ and has rank $d=[R:F]-1$. This can be rephrased as $\varphi|F=0$
and $\Lambda^d(\varphi)=0$ and $\Lambda^{d-1}(\varphi)\neq 0$, a combination of closed and open conditions.

Using universal sheaves on $\prod_{i=1}^s\Grass^{n_i}(R)$, one easily sees that the
combination of these open and closed conditions defines the desired locally closed set $Z$.
\qed

\medskip
In positive characteristics, we have to deal with an additional condition:

\begin{proposition}
\mylabel{second locally closed}
Suppose $p>0$. There is a locally closed set $Z'\subset Z$ such that for each field extension $F\subset K$,
a $K$-valued point $(L_1,\ldots,L_s)$ belongs to $Z'$ if and only if the $K$-algebras $L_i$ admit $p$-bases.
\end{proposition}

\proof
Let $A\subset B$ be a ring extensions that is finite and locally free.
According to \cite{Kiehl; Kunz 1965}, Satz 6 the $A$-algebra $B$ locally admits a $p$-basis if and only if $B^p\subset A$ and 
the $B$-module $\Omega^1_{B/A}$ is locally free.
The former condition means that $g_1^p,\ldots,g_m^p\in A$, where $g_1,\ldots,g_m\in B$ are $A$-algebra generators. 
This in turn means that the $A$-linear map $B^{\oplus m}\ra B/A$ given by 
$(\lambda_i)\mapsto [\sum \lambda_ig_i]$ is zero, which is a closed condition.
Regarding the second condition, the subset inside $\Spec(B)$ where the sheaf  attached to the module of finite presentation 
$\Omega^1_{B/A}$ is locally free is obviously open. Since $\Spec(B)\ra \Spec(A)$ is a closed map,
the subset inside $\Spec(A)$ over which $\Omega^1_{B/A}$ is locally free is likewise open.
 Our assertion follows by applying the preceding observations with the
universal sheaves corresponding to $L_i$.
\qed

\medskip
We remark in passing that a closer analysis reveals that our conditions put a canonical scheme structure
on $Z$ and $Z'$, making them subscheme in the product of Grassmannians $\prod_{i=1}^s\Grass^{n_i}_{R/F}$. 
This, however, plays no role in later applications.

Obviously, the  formation of the locally closed sets $Z$ and $Z'$ commutes with ground field extensions.
By Hilbert's Nullstellensatz, they are uniquely determined by the conditions in the preceding propositions.
Note, however, that they may well be empty.
We now establish    existence results in   two particular cases that are relevant throughout.

\begin{proposition}
\mylabel{existence subrings}
Suppose $R$ admits a $p$-basis of length $r\geq 2$. Then 
there are subalgebras $\Lambda_1,\Lambda_2\subset R$ that admit $p$-bases of length $r-1$,
and satisfy $\Lambda_1\cap \Lambda_2=F$.
\end{proposition}

\proof
Choose a $p$-basis $\omega_1,\ldots,\omega_r\in R$. Then 
$\omega_i'=\omega_i$ for $i\leq r-1$ and $\omega'_r=\omega_r+1$ constitute another $p$-basis,
as one easily sees by taking differentials. It $1+\omega_r$ is not a unit, it must belong
to the maximal ideal of the local Artin ring, and thus $1+\omega'_r$ is a unit.
So without restriction of generality we may assume that $1+\omega_r\in R^\times$.
Now set $s=r-1$  and consider
$$
\Lambda_1=F[\omega_1,\ldots,\omega_s]\quadand \Lambda_2=F[\omega_1+\omega_1\omega_r,\ldots,\omega_s+\omega_s\omega_r].
$$
Obviously the generators form a $p$-basis for the first subring. For the second, we observe 
 $(\Lambda_2)^p\subset F$ and $d(\omega_i+\omega_i\omega_r)=d(\omega_i)(1+\omega_r) + \omega_id(\omega_r)$.
These differentials are linearly independent
in $\Omega^1_{R/F}$, because $1+\omega_r\in R^\times$,  and we conclude with \cite{Kiehl; Kunz 1965}, Satz 6 that  the generators form a $p$-basis.

To see  $\Lambda_1\cap \Lambda_2=F$, suppose we have   $P,Q\in F[T_1,\ldots,T_s]$ with all exponents bounded by $p-1$, such that
$P(\omega_1,\ldots,\omega_s) = Q(\omega_1+\omega_1\omega_r,\ldots,\omega_s+\omega_s\omega_r)$.
Write 
$$
Q(T_1+T_1T_r,\ldots,T_s+T_sT_r) =\sum_{i=0}^s T_r^i Q_i(T_1,\ldots,T_s) 
$$
Then $Q_0=Q$, and the  Taylor expansion (compare \cite{A 4-7}, Chapter IV, \S4, No.\ 5) gives $Q_1=\sum_{i=1}^sT_i\frac{\partial Q}{\partial T_i}$.
Combining these equations and comparing coefficients gives the equation $\sum_{i=1}^s\omega_i\frac{\partial Q}{\partial T_i}(\omega_1,\ldots,\omega_s)=0$
in $R$. Since the $\omega_i$ are $p$-independent and the exponents in $Q$ are bounded by $p-1$ we actually have 
$\sum_{i=1}^sT_i\frac{\partial Q}{\partial T_i}=0$ in $F[T_1,\ldots,T_s]$.
If $Q$ is homogeneous of degree $0\leq d\leq p-1$, the polynomial $\sum_{i=1}^sT_i\frac{\partial Q}{\partial T_i}$ coincides with
$dQ$, by Newton's Formula. For general $Q$, one may apply this fact to the homogeneous summands,
infers that $Q$ is constant, and thus $\Lambda_1\cap\Lambda_2=F$.
\qed

\medskip
The preceding result is a bit paradoxical,  because if $R$ is a field, the   subspaces
 $\Omega^1_{\Lambda_i/F}\otimes_{\Lambda_i}R$ inside the vector space $\Omega^1_{R/F}$ have an intersection 
of dimension at least $2(r-1)-r=r-2$. 
In the above  proof, we actually have not fully used that the last member $\omega_r$ is part of the $p$-basis. This leads to the following:

\begin{proposition}
\mylabel{existence dual subrings}
Suppose   $R$  admits a $p$-basis of length $r$,
and let $0\leq s\leq r$.
Then inside the ring of dual numbers $R[\epsilon]$, there are subalgebras of the form $\Lambda=L+H\epsilon$  such that 
the following holds:
\begin{enumerate}
\item The $F$-algebra $L$ admits a $p$-basis of length $s$.
\item The   residue class map    $L\ra R[\epsilon]/(\epsilon)=R$ is injective.
\item As $L$-module, the quotient $\epsilon R/\epsilon H$ is free of rank one.
\item The intersection $R\cap\Lambda$ inside $R[\epsilon]$ coincides with $F$.
\end{enumerate} 
\end{proposition}

\proof
Choose a $p$-basis $\omega_1,\ldots,\omega_r\in E$ and set $L=F[\omega_1+\omega_1\epsilon,\ldots,\omega_s+\omega_s\epsilon]$.
Then $(\omega_i+\omega_i\epsilon)^p=0$. As in the preceding proof, one verifies that   the generators form a $p$-basis,
and that the  intersection $R\cap L$ inside $R[\epsilon]$ coincides with $F$. The former gives (i), and the latter ensures (ii).
 
To proceed, note that the $F$-module $R\epsilon$ is free, and the $p^r$ monomials 
\begin{equation}
\label{basic monomials}
P_\nu=\prod_{i=1}^r (\omega_i + \omega_i\epsilon)^{\nu_i}\cdot \epsilon = \prod_{i=s+1}^r \omega_i ^{\nu_i}\epsilon\quad  (0\leq \nu_i\leq p-1)
\end{equation}
form   a basis. We see that $R\epsilon$ remains free as an $L$-module, with a basis formed by the monomials $P_\nu$
with     $\nu_1=\ldots=\nu_s=0$.
In particular, the submodule $L\epsilon\subset R\epsilon$ is generated by the monomial with $\nu_1=\ldots=\nu_r=0$.
It  admits a complement, for example the  the submodule $H\epsilon\subset R\epsilon$
generated by the monomials where  $\nu=(\nu_1,\ldots,\nu_r)$ is non-zero with $\nu_1=\ldots=\nu_s=0$.
Clearly, $E\epsilon/H\epsilon = L\epsilon$ is free of rank one, hence (iii).

It remains to verify the statement (iv) about the resulting subring $\Lambda=L+H\epsilon$.
For this we use the  decomposition $R[\epsilon]=R\oplus R\epsilon$ of abelian groups and the corresponding projections.
Note that these respect the module structures over $R$ and in particular over $F$, but not over $L$.
Clearly, $\pr_1|L$ is injective. It follows that $L\cap H\epsilon =0$, hence the sum $\Lambda=L+H\epsilon$ is direct,
and $\dim_F(\Lambda)=p^s + p^{r-s}$.
The intersection $R\cap \Lambda$ is isomorphic to  the kernel of $\pr_2|\Lambda$. To understand the latter,
we examine the image $\pr_2(\Lambda)\subset R\epsilon$, which is generated by the monomials $P_\nu$ of the form
$$
\nu=(\nu_1,\ldots,\nu_s,0,\ldots,0)\quador \nu=(0,\ldots,0,\nu_{s+1},\ldots,\nu_r).
$$
These are $p^s+(p^{r-s}-1)$ members of the $F$-basis \eqref{basic monomials}, which gives the rank of  $\pr_2|\Lambda$.
It follows that 
$
\dim_F(R\cap\Lambda)=\dim_F(\Lambda) - \rank (\pr_2|\Lambda) = 1
$.
In turn, the inclusion $F\subset R\cap\Lambda$ must be an equality.
\qed

\section{Constructions    via  pairs of subfields}
\mylabel{Constructions subfields}

In this section we describe our first construction that leads to the regular  genus-one curves.
In contrast to later constructions, it does not rely on more advanced techniques from commutative algebra,
and has some immediate geometric appeal. The drawback is that it only works in characteristic two.

So let $F$ be an imperfect ground field of characteristic $p=2$,
and $F\subset E$ be a height-one extension of     degree $[E:F]=p^r$ for some $r\geq 1$.
Let $X$ be an regular curve with $H^0(X,\O_X)=E$. Note that we regard this as a scheme over $F$,
such that $h^0(\O_X)=p^r$. In this section we assume
that there are two $E$-valued point $a'\neq a''$ on $X$, and will  describe  certain denormalizations $f:X\ra Y$
that are relevant in the   construction  of regular genus-one curves.

Let $L',L''\subset E$ be   subextensions of  degrees $p^{r-1}$, with $L'\cap L''=F$.
Note that according to Proposition \ref{existence subrings}, such  fields indeed exist. Set 
$$
A=\Spec(E\times E)=\{a',a''\}\quadand B=\Spec(L'\times L'')=\{b',b''\},
$$
and consider the resulting morphism $f:A\ra B$  with $b'=f(a')$ and $b''=f(a'')$. 
In turn, the cocartesian diagram
$$
\begin{CD}
A	@>>>	X\\
@VVV		@VVfV\\
B	@>>> 	Y
\end{CD}
$$
defines a new integral curve $Y$ with singular locus $\Sing(Y)=\{b',b''\}$.  
According to Proposition \ref{fundamental facts}, we have $h^0(\O_Y)=1$ and $h^1(\O_Y)=h^1(\O_X)+1$. Moreover, the local rings of $Y$ are 
geometrically unibranch and Gorenstein, the latter relying on our assumption $p=2$.

The inclusion $L'\cup L''\subset E$ is never an equality (\cite{A 4-7}, Chapter V, \S7, No.\ 4, Lemma 1), so there is some $\omega\in E$ not belonging
to either of the subfields. This $\omega\in E$ is a $p$-basis over both $L'$ and $L''$,
so the scalar $\lambda=\omega^2$ belongs to $L'\cap L''=F$, and is a non-square in both intermediate fields $L'$ and $L''$.

We next examine the complete local noetherian rings at the singularities.
The situation is symmetric, so we fix some $b\in \{b',b''\}$,   write $L\in\{L',L''\}$ for the residue field,
and $a\in\{a',a''\}$   for the corresponding point on $X$.
Fix a uniformizer $u\in \O_{X,a}^\wedge$, such that $\O_{X,a}^\wedge=E[[u]]$.
Note that the subring $\O_{Y,b}^\wedge$ comprises the formal power series $\sum\gamma_iu^i$ whose constant term $\gamma_0$ belongs to the subfield $L\subset E$.

\begin{proposition}
\mylabel{local rings 1,1}
Disregarding the $F$-structure, the complete local noetherian ring $\O_{Y,b}^\wedge$  is   isomorphic to $L[[x,y]]/(x^2-\lambda y^2)$.
\end{proposition}

\proof
Clearly the linear polynomials $\omega u$ and $u$ belong to the subring $\O_{Y,b}^\wedge$, and 
satisfy the relation $(\omega u)^2=\lambda\cdot (u)^2$.  The inclusion  $Lu^0\subset \O_{Y,b}^\wedge$ is a field of representatives, depending on the choice of uniformizer.
The assignment $x\mapsto \omega u$ and $y\mapsto u$ now
defines a homomorphism $\varphi:L[[x,y]]/(x^2-\lambda y^2)\ra \O_{Y,b}^\wedge$.
The ring $R= L[[x,y]]/(x^2-\lambda y^2)$ is a one-dimensional integral domain, the latter because 
 $\lambda y^2\in L[[y]]$ is a  non-square. 

Suppose we have some $\gamma\in E$ and $i\geq 1$.
Write $\gamma=\alpha+\beta\omega$ with $\alpha,\beta\in L$. Then
$\gamma u^i = \alpha u^i + \beta u^{i-1}\cdot \omega u$ belongs to the image of $\varphi$,
and one easily infers that the map is surjective. 
Since $R$ is integral, of the same dimension as $\O_{Y,b}$, we can apply Krull's Principal Ideal Theorem
and infer that $\varphi$ is bijective.
\qed

\medskip
We now specialize to the case that
$$
X=\PP^1_E=\Proj E[T_0,T_1]=\Spec E[u]\cup \Spec E[u^{-1}]
$$
is the projective line over $E$, with inhomogeneous coordinate $u=T_1/T_0$.
Without loss of generality $a'=(1:0)$ and $a''=(0:1)$. The arguments in the preceding proof work with polynomials instead of power series.
It follows that the affine open covering of $X$  induces an affine open covering
$$
Y= \Spec L'[u,\omega u]\cup \Spec L''[u^{-1},\omega u^{-1}],
$$
and the coordinate ring for the overlap  $U=Y\smallsetminus\{b,b'\}$ is the Laurent polynomial ring $E[u^{\pm}]$.
The differential $d\omega\in\Omega^1_{E/F}$ is non-zero, because it is non-zero in the quotients $\Omega^1_{E/L'}$.
Thus we find a basis of the form $d\omega_1,\ldots,d\omega_r\in \Omega^1_{E/F}$, with $\omega_1=\omega$ and some $\omega_i\in E$.
Then 
$du,d\omega_1,\ldots,d\omega_r$ form a basis of $\Omega^1_{U/F}$,
and we write $\partial/\partial u, \partial/\partial\omega_1,\ldots\partial/\partial\omega_r$ for the dual basis,
which lives in the tangent sheaf $\Theta_{U/F}$.
Consider the vector field
$$
D= \partial/\partial u +  \omega u^{-1}\partial/\partial\omega= u^{-2}\partial/\partial u^{-1} +  \omega u^{-1}\partial/\partial\omega,
$$
a priori defined outside $\Sing(Y)=\{b',b''\}$.
 One easily computes  that  the induced derivation     on the function field  has
\begin{equation}
\label{derivation values}
D(u^{-1})=u^{-2},\quad D(u)=1,\quad D(\omega)=\omega u^{-1},\quad D(\omega u^{-1}) = D(\omega u) =0.
\end{equation}
Consequently, $D:E(u)\ra E(u)$ stabilizes  the two local rings $\O_{Y,b'}$ and $\O_{Y,b''}$ inside $\Gamma(U,\O_Y)$,
and   $D^{[2]}=D\circ D$ is the zero map. It thus defines a non-zero additive $p$-closed global vector field on $Y$,
which corresponds to an inclusion $\alpha_p\subset\Aut_{Y/F}$, in other words, a faithful $\alpha_p$-action.

Let $Z=\alpha_p\cdot Z_0$ be the orbit of the reduced closed subscheme $Z_0=\{b',b''\}$. 
The geometry of our $\alpha_p$-action   has completely different features  at the two singular points:

\begin{proposition}
\mylabel{action 1,1}
In the above setting, the $\alpha_p$-scheme $Z$ is equivariantly isomorphic to the disjoint union $(\alpha_p\times\Spec L')\cup (\Spec L'') $,
where the action   is described on the corresponding functors of $k$-algebras by  the formula
$$
\sigma\cdot (\tau,z') = (\sigma\tau,z')\quadand\sigma\cdot z''=z''.
$$
Moreover, the Weil divisor $Z\subset Y$ is Cartier at the point $b'$, but not at $b''$.
\end{proposition}

\proof
Let $\shI_0\subset \O_Y$ be the sheaf of ideals for $Z_0$. Then the kernel $\shI$
of the additive map $\shI_0\stackrel{D}{\ra}\O_{Z_0}$ corresponds to the scheme of orbits $Z\subset Y$.
From \eqref{derivation values} we see that 
 $\maxid_{b''}$ is $D$-stable. It follows that the inclusion $Z_0\subset Z$ is an equality at $b''$.
Since this point is singular, the inclusion $Z\subset Y$ is not Cartier at $b'$.
Furthermore, one sees that the induced derivation on the residue field $L''=\kappa(b'')$ is trivial.

On the other hand,  $\maxid_{b'}$ is not $D$-stable, because $D(u)=1$.
Furthermore we have $\omega u\in \shI_{b'}$,
and see  $\O_{Y,b'}/(\omega u) = L'[x,y]/(x^2-\lambda y^2, x)= L'[y]/(y^2)$. Using that   $\maxid_{b'}/(\omega u)$ is a simple module
we infer $\shI_{b'}=(\omega u)$. This shows that $Z\subset Y$ is Cartier near $b$.
Furthermore, \eqref{derivation values} reveals  that the  $\alpha_p$-action on $Z$ at $b$ is given  by the derivation $\partial/\partial y$
on the coordinate ring $L'[y]/(y^2)=F[y]/(y^2)\otimes L'$, and the assertion follows.
\qed

\medskip
By symmetry, the above analysis also applies to the vector field
$$
\tilde{D}=\partial/\partial u^{-1} + \omega u\partial/\partial\omega= u^2\partial/\partial u  +  \omega u\partial/\partial\omega,
$$
with the roles of $b',b''$ interchanged. One easily checks that  $[D,\tilde{D}]=0$, and that the two derivations are $F$-linearly independent.
In turn, the two global vector fields $D,\tilde{D}$ define a two-dimensional restricted Lie algebra $\lieg=k^2$ inside $H^0(Y,\Theta_{Y/F})$,
with trivial bracket and $p$-map.
In other words, we have an action of $G=\alpha_p\times\alpha_p$ on $Y$.

Now $  F',F''$  be two  simple height-one extensions of $F$, and 
endow their spectra  with the   structure of an $\alpha_p$-torsor.  
The affine scheme $T=\Spec(F'\otimes F'')$ becomes a torsor with respect to the  group scheme $G=\alpha_p\times\alpha_p$.
In turn, we obtain the twisted form
$$
\tilde{Y}=T\wedge^G Y = G\backslash (T\times Y).
$$
This defines another genus-one curve $\tilde{Y}$. Now recall that   any purely inseparable field extensions of $F$
uniquely embed into $F^\alg$, so the concept of \emph{linear disjointness} applies
(\cite{A 4-7}, Chapter V, \S 2, No.\ 5), even without giving an ambient field.

\begin{proposition}
\mylabel{regular 1,1}
Suppose   the three   height-one extensions $F\subset E,F',F''$ are linearly disjoint. 
Then the   genus-one curve $\tilde{Y}$ is regular. Moreover, its Frobenius base-change $\tilde{Y}^{(p)}$ is 
is isomorphic  to some standard model $C^{(1,1)}_{r,F,\Lambda}$.
\end{proposition}

\proof
The orbit $Z=G\cdot Z_0$ of the reduced scheme $Z_0=\{b',b''\}$ is $G$-stable.
It is 
equivariantly isomorphic to the disjoint union $(\alpha_p\times\Spec L')\cup (\alpha_p\times \Spec L'')$,
according to Proposition \ref{action 1,1}, and  the $G$-action is given by
$$
(\sigma_1,\sigma_2)\cdot (\tau,z') = (\sigma_1\tau,z')\quadand (\sigma_1,\sigma_2)\cdot (\tau,z'') = (\sigma_2\tau,z'').
$$
In turn, the twisted form $\tilde{Z}=Z\wedge^GT$ has coordinate ring $(F'\otimes L')\times (F''\otimes L'')$.
This Artin ring is regular, because by assumption the tensor products are fields.
According to  Proposition \ref{action 1,1}, the inclusion $Z\subset Y$ is an effective Cartier divisor, so the same holds
for $\tilde{Z}\subset\tilde{Y}$. In turn, the scheme $\tilde{Y}$ is regular at the two points of the divisor.

Finally, consider the open set $U=\Reg(Y)$. The ensuing twisted form $\tilde{U}$ comes  with
a faithfully flat morphism from the spectrum of $E[u^{\pm 1}]\otimes F'\otimes F''$,
which can be seen as the Laurent polynomial ring over $E\otimes F'\otimes F''$.
Again, our assumption ensures that this ring is regular, hence $\tilde{U}$ is regular.
Thus Lemma \ref{twisting removes singularities} ensures that the genus-one curve $\tilde{Y}$ is regular.

The Frobenius pull-back $E\otimes_FF$ is isomorphic to $F[W_1,\ldots,W_r]/(W_1^p,\ldots,W_r^p)$,
and likewise for the field extensions $L',L'',F',F''$. From this we see $Y^{(p)}$ is isomorphic
to some standard model $C^{(1,1)}_{r,F,\Lambda}$.
Moreover,  $T^{(p)}$ becomes the trivial torsor
with respect to $G=G^{(p)}$, thus $Y$ and its twisted form $\tilde{Y}$ become isomorphic after
Frobenius pull-back. Summing up, $\tilde{Y}^{(p)}\simeq C^{(1,1)}_{r,F,\Lambda}$.
\qed


\section{Constructions with fields of representatives}
\mylabel{Constructions fields representatives}

In this section $F$ is an imperfect ground field of arbitrary characteristic $p>0$.
Again  $X$ denotes  a regular curve where $E=H^0(X,\O_X)$ is a height-one extension of degree $[E:F]=p^r$ for some $r\geq 1$.
But now we fix only one $E$-valued point $a\in X$, and will describe another denormalization
 that is relevant in the construction of  regular genus-one curves.
Consider the local Artin ring $\O_{X,a}/\maxid_a^2$. Throughout we will identify $E=\Gamma(\O_X)$
with its image in $\O_{X,a}/\maxid_a^2$, which  is the  \emph{standard field of representatives}.

We now fix \emph{another field of representatives}   $L\subset \O_{X,a}/\maxid_a^2$, having the property
$\Gamma(\O_X)\cap L=F$.  According to Proposition \ref{existence dual subrings} such fields indeed exist.  Note that
the projection to the residue field gives a canonical identification between $E=\Gamma(\O_X)$ and $L$, but
it is crucial to view them  as \emph{different} fields inside $\O_{X,a}/\maxid_a^2$. Set 
$$
A=\Spec(\O_{X,a}/\maxid_a^2)=\{a\}\quadand B=\Spec(L)= \{b\}.
$$
The ensuing cocartesian diagram
\begin{equation}
\label{cocartesian first}
\begin{CD}
A	@>>>	X\\
@VVV		@VVfV\\
B	@>>> 	Y.
\end{CD}
\end{equation}
defines a new integral curve $Y$, with  $\Sing(Y)=\{b\}$. 
Arguing as in Proposition \ref{standard models}, one sees that 
$h^0(\O_Y)=1$ and  $h^1(\O_Y)=h^1(\O_X)+1$. Furthermore, $Y$ is locally unibranch and Gorenstein.

We now seek to understand the complete local ring $\O^\wedge_{Y,b}$ at the $E$-valued singularity $b\in Y$.
To this end we fix a uniformizer $u\in \O_{X,a}^\wedge$, such that $\O_{X,a}^\wedge=E[[u]]$ as $F$-algebra.
This turns $\O_{X,a}/\maxid_a^2= E[\epsilon]$
into a ring of dual numbers, where $\epsilon $ is the   class of the uniformizer.
Choose  a  $p$-basis $\omega_1,\ldots,\omega_r\in E$.
The   residue class bijection $L\ra\kappa(b)= \kappa(a)=E$ takes the form
$$
\omega_i+\alpha_i\epsilon\longmapsto \omega_i 
$$
for some uniquely determined $\alpha_i\in E$, and the elements $\omega_i+\alpha_i\epsilon\in E[\epsilon]$ form a $p$-basis for the subfield $L$. 
The cartesian square 
\begin{equation}
\label{ring representatives}
\begin{CD}
E[\epsilon]						@<<<	E[u]\\
@AAA							@AAA\\
F[\omega_1+\alpha_1\epsilon,\ldots,\omega_r+\alpha_r\epsilon]	@<<<	R
\end{CD}
\end{equation}
defines a subring $R\subset E[u]$, which comprises the polynomials  $\beta_0+\beta_1u+\ldots+\beta_nu^n$ whose  truncation $\beta_0 + \beta_1\epsilon$
is  a polynomial     in  the  $\omega_i+\alpha_i\epsilon$. In light of the Taylor expansion, this precisely means 
\begin{equation}
\label{taylor expension}
\beta_0=P(\omega_1,\ldots,\omega_r)\quadand \beta_1= \sum_{i=1}^r \alpha_i\frac{\partial P}{\partial T_i}(\omega_1,\ldots,\omega_r)
\end{equation}
for some   $P\in F[T_1,\ldots,T_r]$.  

Clearly the $R$-algebra $E[u]$ is a finite, hence 
$\maxid=R\cap uE[u]$ is a maximal ideal, having  residue field   $ R/\maxid=L$. 
Moreover, 
$R$ is a finitely generated $F$-algebra.
It appears that the number of  generators is excessive,
so to improve the situation we have to pass to the formal completion
$\widehat{R}=\invlim_n R/\maxid^n$. The latter sits in a cartesian diagram analogous to \eqref{ring representatives},
with $E[[u]]=\O_{Y,a}^\wedge$ in the top right corner,  
and is canonically identified with $\O_{Y,b}^\wedge$. 

\begin{proposition}
\mylabel{local ring representative}
Disregarding the $F$-structure, the complete local ring $\widehat{R}=\O_{Y,b}^\wedge$ is   isomorphic to $E[[x,y]]/(x^3-y^2)$.
\end{proposition}

\proof 
Choose a section $s$ for the residue class projection $\widehat{R}\ra L$. It takes the form $s(\omega_i) = \omega_i+\alpha_iu + \ldots$.
Set $\widehat{R}_0=L[[x,y]]/(x^3-y^2)$. The relation is irreducible   because $x^3$ is not 
a square in   $L[[x]]$.
Hence $\widehat{R}_0$ is a one-dimensional integral domain.
Clearly,     $u^2,u^3\in E[[u]]$ belong to $\maxid=R\cap uE[u]$
and satisfy $(u^2)^3=(u^3)^2$.  Together with our field of representatives $s(L)\subset\widehat{R}$,
the assignment $x\mapsto u^2$ and $y\mapsto u^3$ defines  a map 
$\varphi:\widehat{R}_0\ra \widehat{R}$. 

We claim that each formal power series  $\sum_{i\geq 2}\lambda_i u^i$ with coefficients from $E$ belongs to $\Image(\varphi)$.
Write $\lambda_i=P_i(\omega_1,\ldots,\omega_r)$ as a $p$-truncated polynomial with coefficients from $F$.
Then 
$$
\lambda_iu^i = P_i(\omega_1,\ldots,\omega_r)u^i=P_i(s(\omega_1),\ldots,s(\omega_r))u^i + \ldots,
$$ 
where the missing terms have order at least $ i+1$.
The   factor  $P_i(s(\omega_1),\ldots,s(\omega_r))$ belongs to $\varphi(s(L))$,
whereas $u^i=\varphi(y^{i/2})$ or $u^i=\varphi(xy^{(i-3)/2})$, depending on the parity of $i$. In any case, the exponent  in the $y$-term
grows linearly with  $i\geq 2$, hence successive substitutions   reveal that  $\sum_{i\geq 2}\lambda_i u^i$
belongs to $\Image(\varphi)$.

In light of the cartesian diagram \eqref{ring representatives}, the maximal ideal of $\widehat{R}$ consists of these $\sum_{i\geq 2}\lambda_i u^i$,
and it follows that the map $\varphi:\widehat{R}_0\ra \widehat{R}$ is surjective. 
The rings have the same dimension and $\widehat{R}_0$ is integral. Using Krull's Principal Ideal Theorem, we infer
that $\varphi$ is bijective.
\qed

\medskip
We next seek to understand  $\widehat{R}=\O_{Y,b}^\wedge$ as   $F$-algebra. Clearly, the $r+2$ polynomials
\begin{equation}
\label{generators}
u^3\quadand  u^2\quadand \omega_i+\alpha_iu \quad (1\leq i\leq r) 
\end{equation}
are contained in $R$. The  $p$-powers  $\lambda_i=\omega_i^p$ and $\mu_i=\alpha_i^p$ belong to the ground field $F$, and the 
 above generators satisfy the $r+1$ \emph{obvious relations}
\begin{equation}
\label{relations}
(u^3)^2-(u^2)^3=0\quadand (\omega_i+\alpha_iu)^p- \lambda_i-\mu_iu^p=0 \quad (1\leq i\leq r).
\end{equation}
Here one has to  rewrite, for $p\geq 5$, the factor $u^p$ as $(u^3)\cdot (u^2)^{(p-3)/2}$, to get an expression   in terms of the generators.
By abuse of notation, we also regard the generators in \eqref{generators} as indeterminates, and write $\tilde{R}$ for the resulting
polynomial ring, formally completed with respect to the ideal $(u^2, u^3)$.
From universal properties  we get a continuous homomorphism
$$
\varphi:\tilde{R}=F[\omega_1+\alpha_1u,\ldots,\omega_r+\alpha_ru][[ u^2,u^3]]\lra \widehat{R}.
$$
The following gives the desired description  as completed $F$-algebra:
 
\begin{proposition}
\mylabel{affine model representatives}
The above map is surjective,   the  obvious relations in \eqref{relations}
form a regular sequence in   $\tilde{R}$, and they generate the ideal  $\idealp=\Kernel(\varphi)$ of all relations. 
\end{proposition}

\proof
Let $\ideala\subset \tilde{R}$ be the ideal generated by the obvious relations in \eqref{relations}.
The main task is to verify that the  inclusion $\ideala\subset\idealp$ is an equality.
One easily checks  
$$
\tilde{R}/(\ideala+\tilde{R}u^2)=L[u^3]/(u^6)\quadand (\tilde{R}/\ideala)_{u^2}=L[[u]]_{u^2} = L((u)).
$$
Hence there are exactly two primes $\primid_1,\primid_2\subset\tilde{R}$ containing $\ideala$,
one being the maximal ideal $\primid_1=\ideala + \tilde{R}u^2+ \tilde{R}u^3$, which has height $r+2$.
By Krull's Principal Ideal Theorem, every minimal prime containing
$\ideala$ has height at most $r+1$. 
We infer $\dim(\tilde{R}/\ideala) = 1$, which coincides with the difference $ (r+2)-(r+1)$. It follows that
 the $r+1$ obvious relations  from    \eqref{relations} in the $r+2$ generators in \eqref{generators} form a regular sequence, and  
$\tilde{R}/\ideala$ is Cohen--Macaulay (\cite{SP}, Tag 02JN).
In our situation, the latter simply means that the quotient has no embedded primes. 
Clearly $R_{u^2}=E[u^{\pm 1}]$, hence  $\ideala\subset\primid$ becomes an equality after localization of $u^2$.
Thus $\primid/\ideala$ is an $\tilde{R}$-module of finite length. It must vanish, because
it  coincides with the kernel of $\tilde{R}/\ideala\ra R$, and $\tilde{R}/\ideala$ has no embedded primes.

It remains to check that the inclusion $\tilde{R}/\primid\subset R$ is an equality.
One easily sees that it becomes an equality after localization of $u^2$, and also after formally completing
with respect to $(u^2)$, and the result follows.
\qed

\medskip
Recall that the \emph{Fitting ideals} $\Fitt_i(M)$ for an $A$-module $M$ having a presentation $A^n\ra A^m\ra M\ra 0$
are generated by the $(m-i)$-minors of the presentation matrix, compare the discussion in \cite{Brion; Schroeer 2022}, Section 3.
The above analysis enables us to understand  K\"ahler differentials and their Fitting ideals,
a method already used in loc.\ cit.\ in connection with automorphism group schemes.

\begin{proposition}
\mylabel{kaehler representatives}
The coherent sheaf $\Omega^1_{Y/F}$ has rank $r+1$, and the following are equivalent:
\begin{enumerate}
\item The torsion-free sheaf $\Fitt_{r+1}(\Omega^1_{Y/F})$  is locally free.
\item The characteristic satisfies $p\leq 3$.
\end{enumerate}
If these equivalent conditions hold, $\Omega^1_{Y/F}$ is locally free near $b\in Y$.
\end{proposition}

\proof
First note that  on the regular locus $Y\smallsetminus\{b\}$, any torsion-free sheaf is locally free.
So we only have to understand the situation for the local ring $\O_{Y,b}$
or its formal completion $\O_{Y,b}^\wedge$.  We get     identifications
$$
\Omega^1_{Y/F,b}\otimes_{\O_{Y,b}}\O_{Y,b}^\wedge=
(\Omega^1_{Y/F,b})^\wedge_{\maxid_b} = 
(\Omega^1_{R/F})_\maxid^\wedge = 
\Omega^1_{R/F}\otimes_{R}\widehat{R},
$$
the outer ones by \cite{AC 1-4}, Chapter III, \S3, No.\ 4, Theorem 3, the middle one  by
the compatible equalities  $\O_{Y,b}/\maxid_b^{n+1}=R/\maxid^{n+1}$ stemming from   Proposition \ref{affine model representatives}. According to
Lemma \ref{presentation  kaehler differentials} below, the finitely presented module $(\Omega^1_{R/F})^\wedge_\maxid$ 
has a presentation with   generators $d(u^2),d(u^3),d(\omega_i+\alpha_iu)$,
and relations obtained by expressing the  differentials of the   relations
\eqref{relations}  in terms of the generators.
Hence   $\Omega^1_{R/F}\otimes_R\widehat{R}$ has a presentation matrix in block form
$$
P=
\left(
\begin{array}{c|ccc}
2u^3 	& P_1 	& \cdots & P_r \\
-3(u^2)^2 	& Q_1 	& \cdots & Q_r \\
\cline{1-4}
0	& 	&0\\
\end{array}
\right)
\in\Mat_{(r+2)\times(r+1)}(\widehat{R}).
$$
For $p\geq 5$, the entries in the upper right block are the coefficients from  the differential of the relation
$(\omega_i+\alpha_iu)^p-\lambda_i- \mu_i\cdot (u^3)\cdot (u^2)^{(p-3)/2}$, so
$$
P_i=-\mu_i(u^2)^{(p-3)/2}\quadand Q_i=-\frac{p-3}{2}\mu_i(u^3)  (u^2)^{(p-5)/2}.
$$
For  the remaining     characteristics we get 
$$
(P_i,Q_i)=\begin{cases}
(0,-\mu_i)	& \text{if $p=2$;}\\
(-\mu_i,0)	& \text{if $p=3$.}
\end{cases}
$$
One quickly computes that in all cases the 2-minors $P_iQ_j-Q_iP_j$ and
$2u^3Q_i+3(u^2)^2P_i$ vanish, and immediately sees that  the presentation matrix has non-zero entries.
For the Fitting ideals, this means  $\Fitt_r(\Omega^1_{Y/F})_b=0$ and $\Fitt_{r+1}(\Omega^1_{Y/F})_b\neq 0$, and it follows that 
the sheaf $\Omega^1_{Y/F}$ has rank $r$.

Suppose now $p=3$. Then all but the first row of $P$ vanish, and this row is
$-(u^3, \mu_1,\ldots,\mu_r)$. Our standing assumption $\Gamma(\O_X)\cap L=F$ ensures
that the $\mu_i=\alpha_i^p$ are non-zero, hence $\Fitt_{r+1}(\Omega^1_{Y/F})_b$ is the unit ideal,
and in particular free as a module.
Thus the numerical function $y\mapsto \dim_{\kappa(y)} (\Omega^1_{Y/F}\otimes\kappa(y))$ is constant on $\{b,\eta\}$,
and it follows that the stalk $\Omega^1_{Y/F,b}$ is free.
The argument for $p=2$ is analogous. This establishes (ii)$\Rightarrow$(i).

It remains to prove (i)$\Rightarrow$(ii).
Recall that    $\idealb=\Fitt_{r+1}(\Omega^1_{R/F} \otimes_R\widehat{R})$ is generated by the entries of
the presentation matrix $P$, and write  $\idealb'$ for the induced ideal in the normalization $\widehat{R}'=E[[u]]$.
Let us first consider the case  $p\geq 7$, which ensures $2,3\in F$ are non-zero, and $(u^2)^2\mid P_i$ and $u^3\mid Q_i$.
Then  $\widehat{R}'/\idealb'=E[u]/(u^3)$, and the relation \eqref{relations} reveals $\widehat{R}/\idealb= L[u^2]/(u^4)$.
For  $p=5$ one  similarly has  $\idealb=(u^3,u^2)$, giving $\widehat{R}'/\idealb'=E[u]/(u^2)$ and 
$R/\idealb= L$. In both cases 
\begin{equation}
\label{degree jumps}
\dim_F(\widehat{R}/\idealb) <\dim_F(\widehat{R}'/\idealb').
\end{equation}
Seeking a contradiction, we now assume that  that $p\geq 5$, and that  $\idealb$ is locally free. 
Being a non-zero ideal, it must be invertible, and thus defines an effective Cartier divisor $D\subset Y$
supported at the singularity $b\in Y$. By the Degree Formula, it has the same degree
as its preimage $f^{-1}(D)\subset X$, in contradiction to \eqref{degree jumps}.
\qed

\medskip
From now on, we assume $p\leq 3$, and  specialize to the case that 
$$
X=\PP^1_E=\Proj E[T_0,T_1]=\Spec E[u]\cup \Spec E[u^{-1}]
$$
is the projective line over our height-one extension $F\subset E$, with inhomogeneous coordinate $u=T_1/T_0$,
and that our chosen point is the origin  $a=(0:1)$. This gives an affine open covering
$$
Y=\Spec R \cup \Spec E[u^{-1}].
$$

Suppose now that there is $D\in H^0(Y,\Theta_{Y/F})$ with either $D^{[p]}=0$ or $D^{[p]}=D$,   that the $\omega_k+\alpha_k u$ belong to $\Kernel(D)$,
and that $D$ does not stabilized the maximal ideal $\maxid_b\subset\O_{Y,b}$.
According to Corollaries \ref{vector field p=2} and \ref{vector field p=3} below, such   global vector fields indeed exists provided that the
field of representatives $L$ is chooses in a special way.
Our  $D\in H^0(Y,\Theta_{Y/F})$ corresponds to a faithful action of $G=\alpha_p$ or $G=\mu_p$.
Let $F\subset F'$ be a simple height-one extension, endow its spectrum $T=\Spec(F')$ with the structure of an $G$-torsor,
and consider the resulting twisted form $\tilde{Y}=T\wedge^G Y$.
 
\begin{proposition}
\mylabel{regular representative}
Assumptions as above. If the height-one extensions $ E,F'$ are linearly disjoint over $F$, then the genus-one curve $\tilde{Y}$ is regular.
Moreover, its Frobenius base-change $\tilde{Y}^{(p)}$ is isomorphic to some standard model $C^{(1)}_{r,F,\Lambda}$.
\end{proposition}

\proof
Consider the singular locus $Z_0=\{b\}$ as reduced closed subscheme, and write $Z=G\cdot Z_0$ for its orbit.
We claim that $Z$ is equivariantly isomorphic to  $G\times\Spec(L)$, where the action is given by
$\sigma\cdot (\tau,z)=(\sigma\tau,z)$, and that the inclusion $Z\subset Y$ is an effective Cartier divisor.
To see this, write $\shI\subset \shI_0$ for the respective sheaf of ideals for the closed subschemes $Z\supset Z_0$.
In order to give a uniform treatment use the Kronecker delta $\delta=\delta_{p,2}$.
Obviously we have
$$
u^p, \omega_1+\alpha_1u,\ldots, \omega_r+\alpha_ru\in \Kernel(D)\quadand D(u^{2+\delta}) = \pm 1.
$$
It follows that $Z\subset Y$ is defined by  $u^p\in\O_{Y,b}$,
with resulting coordinate ring $L[u^{2+\delta}]/(u^{p(2+\delta)})$.
The induced derivation on the coordinate ring vanishes on the   $\omega_i+\alpha_iu$,
and takes the unit value on the generator $u^{2+\delta}$. It follows that $Z$ is equivariantly isomorphic
to $G\times\Spec(L)$. Thus its twisted form $\tilde{Z}$ has coordinate ring
$F'\otimes L$. We also see that the inclusion $Z\subset Y$ is an effective Cartier divisor.

Next, we consider the complementary open set $U=Y\smallsetminus \{b\}$. Its coordinate ring is the
polynomial ring $E[u^{-1}]$, and its base change becomes $(F'\otimes E)[u^{-1}]$.
Our assumptions ensure that the rings $F'\otimes E$ and hence also $F'\otimes L$ are fields.
Thus Lemma \ref{twisting removes singularities} applies, and we conclude that the curve $\tilde{Y}$ is regular.
One argues as in Proposition \ref{regular 1,1} to see that the Frobenius pull-back $Y^{(p)}$ 
is isomorphic to some standard model $ C^{(1)}_{r,F,\Lambda}$.
\qed

\medskip
In the proof for Proposition \ref{kaehler representatives} we have used a useful  general fact on K\"ahler differentials:
Let $R$ be a finitely generated $F$-algebra, $\maxid$ be a maximal ideal, and $\widehat{R}$ be the resulting formal completion.
Suppose we have elements $f_1,\ldots,f_m\in R$ and $g_1,\ldots,g_n\in\maxid$ such that
the former generate the the residue field $\kappa=R/\maxid$ over $F$, and the latter generate the cotangent space $\maxid/\maxid^2$
over $\kappa$. We then have a homomorphism 
\begin{equation}
\label{before completion}
\varphi:F[x_1,\ldots,x_n,y_1,\ldots,y_n]\lra R
\end{equation}
given by the assignments $x_i\mapsto f_i$ and $y_j\mapsto g_j$.
Now suppose we have polynomials  $h_1,\ldots,h_r$ in the indeterminates $x_i$ and $y_j$  that  generate  the kernel for the induced continuous map
$\widehat{\varphi}:F[x_1,\ldots,x_m][[y_1,\ldots,y_n]]\ra \widehat{R}$.

\begin{lemma}
\mylabel{presentation  kaehler differentials}
In the above situation, the $\widehat{R}$-module $(\Omega^1_{R/F})^\wedge_\maxid$ has a presentation
where the generators are the differentials $df_i,dg_j$  and the relations arise
from writing the $dh_k$ in terms of the generators.
\end{lemma}

\proof
To better conform with the cited literature, we  temporarily  change notation  and set 
$$
A=F\quadand B=F[x_1,\ldots,x_m,y_1,\ldots,y_n]\quadand C=R.
$$
Consider the formal completions
$\widehat{B}=\invlim_i B/\idealb^i$ and $ \widehat{C}=\invlim_i C/\maxid^i$, 
where $\idealb=(y_1,\ldots,y_n)$, and $\maxid\subset C$ is the given maximal ideal.
Our homomorphism \eqref{before completion} becomes   $\varphi:B\ra C$. 
The induced continuous map $\widehat{\varphi}:\widehat{B}\ra\widehat{C}$ is surjective. To see this
one argues as in the local case (\cite{AC 8-9}, Chapter IX, \S2, No.\ 5, Lemma 3),
by observing that the induced map   on associated graded rings is surjective, and applying
\cite{AC 1-4},  Chapter III, \S 2, No.\ 8, Corollary  2 for Theorem 1.

Sept $I=\Kernel(\widehat{\varphi})$.
According to  \cite{EGA IVa}, Chapter 0, Theorem 20.5.12 we have an exact sequence
$$
I/I^2\lra \Omega^1_{\widehat{B}/A}\otimes_{\widehat{B}}\widehat{C}\lra \Omega^1_{\widehat{C}/A}\lra 0,
$$
where the map on the left is given by $x\mapsto dx\otimes 1$. But note  that the modules on the right are usually
not finitely generated. To remedy this, we pass to  $\widehat{\Omega}^1_{\widehat{B}/A} = \invlim \Omega^1_{B_i/A}$ where $B_i=B/\idealb^i$.
This module is separated and complete, the ensuing derivation $d:\widehat{B}\ra\widehat{\Omega}^1_{\widehat{B}/A}$
is continuous, and actually universal for continuous derivation to separated and complete modules.
Likewise, we form  $\widehat{\Omega}^1_{\widehat{C}/A} = \invlim \Omega^1_{C_i/A}$ with  $C_i=C/\maxid^i$.  
As discussed in \cite{EGA IVa}, 20.7.14 we still have a sequence
\begin{equation}
\label{differential sequence}
I/I^2\lra \widehat{\Omega}^1_{\widehat{B}/A}\otimes_{\widehat{B}}\widehat{C}\lra \widehat{\Omega}^1_{\widehat{C}/A}\lra 0.
\end{equation}
Such sequences can be defined for every metrisable topological ring $\widehat{B}$ with a closed ideal $I$.
As  remarked in loc.\ cit.\ 20.7.17 and 20.7.20 the images are still dense in the kernels, but 
otherwise the  above sequence may loose its exactness property.
 We now argue that this does not happen here, by using the work of Kunz:

According to \cite{Kunz 1986}, the ring  $\widehat{B}$ admits a derivation $\widehat{B}\ra\widetilde{\Omega}^1_{\widetilde{B}/A}$
to a finitely generated module that is universal foe derivations to finitely generated modules,
and the same holds for  $\widehat{C}$. Moreover, in our situation  we have canonical identifications
$$
\widetilde{\Omega}^1_{\widetilde{B}/A}=\widehat{\Omega}^1_{\widehat{B}/A}\quadand
\widetilde{\Omega}^1_{\widetilde{C}/A}=\widehat{\Omega}^1_{\widehat{C}/A},
$$
according to loc.\ cit.\ Corollary 12.5. Using loc.\ cit.\ Corollary 11.10 we see that \eqref{differential sequence} is indeed exact.
 
One directly checks that the $\widehat{B}$-modules $\widehat{\Omega}^1_{\widehat{B}/A}$ is freely generated by the $dx_i$ and $dy_j$.
By assumption we have  $I=(h_1,\ldots,h_r)$.
The only remaining task is  to identify $\widehat{\Omega}^1_{\widehat{C}/A}=\widehat{\Omega}^1_{\widehat{R}/F}$
with $(\Omega^1_{R/F})^\wedge_\maxid$. To see this set $R_i=R/\maxid^i$, and consider the exact sequences
$
\maxid^i/\maxid^{2i}\ra \Omega^1_{R/F}/\maxid^i\Omega^1_{R/F} \ra \Omega^1_{R_i/F}\ra 0
$.
The image of the map on the left is denoted by  $M_i$. Passing to inverse limits, we get an exact sequence
$$
\invlim M_i \lra (\Omega^1_{R/F})^\wedge_\maxid \lra \widehat{\Omega}^1_{\widehat{R}/F}\stackrel{\partial}{\lra} R^1\invlim M_i.
$$
Clearly, the transition maps $M_{2i}\ra M_i$ are zero, hence the term on the left vanishes. 
Since the   $M_i$ are modules of  finite length, the inverse system $(M_i)_{i\geq 0}$ automatically satisfy the Mittag-Leffler Condition,
so the connecting map $\partial$ is zero  (\cite{EGA IIIa}, Proposition 13.2.2, confer also \cite{Roos 2006} and \cite{Neeman 2002}).
\qed

\section{Analysis  of the Lie algebra}
\mylabel{Lie algebra}

We keep the setting of the previous section, with the genus-one curve 
$$
Y=\Spec R \cup \Spec E[u^{-1}]
$$
obtained as denormalizations of $X=\PP^1_E$ by introducing a singularity $b\in Y$ on the first chart, for the moment with $p>0$ arbitrary.
We now seek to gather information on   $\lieg=H^0(Y,\Theta_{Y/F})$. This is  a finite-dimensional restricted Lie algebra
contained in  $\Der_F(E[u^{-1}])$,  via the restriction map to the regular locus $U=\Spec E[u^{-1}]$.
Each element  of $\Gamma(U,\Theta_{Y/F})= \Der_F(E[u^{-1}])$ takes the form
$$
D = -P\frac{\partial}{\partial u^{-1}} + \sum_{i=1}^r Q_i \frac{\partial}{\partial\omega_i}=u^2P\frac{\partial}{\partial u} + \sum_{i=1}^r Q_i\frac{\partial}{\partial\omega_i},
$$
with    unique polynomials $P,Q_i\in E[u^{-1}]$. Every such $D$ induces an $F$-derivation of the ring of formal Laurent series $E((u))$,
and we have $D\in\lieg$ if and only if this induced map stabilizes the subring $\widehat{R}=\O_{Y,b}^\wedge$.
We now seek to characterize the condition $D\in\lieg$ in terms of the coefficients of $P$ and $Q_i$.
Let us first corroborate that $\lieg$ is finite-dimensional:

\begin{proposition}
\mylabel{bounds on degrees}
If $D\in\lieg$ then  $\deg(P)\leq 4$ and $\deg(Q_i)\leq p$. In odd characteristics   we actually have $\deg(P)\leq 3$.
\end{proposition}

\proof
One immediately sees 
$$
D(u^2)=2u^3P\quadand D(u^3)=3u^4P \quadand D(\omega_ku^p) = u^pQ_k\quad (1\leq k\leq r).
$$
These      Laurent polynomials  belong to the singular subring $\O_{Y,b}^\wedge$, because this holds for the arguments $u^2,u^3,\omega_ku^p$,
and in particular   lie in $E((u))$. As a consequence, they belong to 
 $E[u]=E[u^{\pm 1}]\cap E((u))$,   and the assertion follows.
\qed

\medskip
Recall that $\omega_1,\ldots,\omega_r\in E$ form a    $p$-basis, hence each $\varphi\in E$ can be uniquely written
as $p$-truncated  polynomial $\varphi=\varphi(\omega_1,\ldots,\omega_r)$ with coefficients from $F$.
In turn, we may regard $E$ as a \emph{$D$-module}, that is, 
a module over the \emph{Weyl algebra} $A_r(F)$ formed with the symbols $x_i,\partial_i$.
Explicitly, the symbols act via 
$x_i\cdot \varphi=\omega_i\varphi$ and  $\partial_i\cdot\varphi = \frac{\partial_i\varphi}{\partial x_i} (\omega_1,\ldots,\omega_r)$.
The $D$-module is annihilated by the central elements $x_i^p-\omega_i^p$, so we may as well work with  the residue class ring
\begin{equation}
\label{quotient weyl algebra}
A_r(F)/(x_1^p-\omega_1^p,\ldots,x_r^p-\omega_r^p) = E[\partial_1,\ldots,\partial_r].
\end{equation}
Note that this is the associative $E$-algebra given by the relations
$\partial_i\omega_j - \omega_j\partial_i =\delta_{ij}$. 

Suppose now $p=2$, and consider derivations $D$ whose coefficients take the form
$$
P(u^{-1})=\lambda_0u^{-4} + \lambda_1u^{-3}+\ldots + \lambda_4\quadand Q_i(u^{-1}) = \mu^{(i)}_0 u^{-2} +  \mu^{(i)}_1 u^{-1} + \mu^{(i)}_2.
$$
Consider, for $1\leq k\leq r$, the  following elements from the field $E$:
$$
\Phi_k = \alpha_k\lambda_2 + \mu_2^{(k)} + \sum_{i=1}^r \mu_1^{(i)} \partial_i\alpha_k \quadand 
\Psi_k = \alpha_k\lambda_3 + \sum_{i=1}^r\mu^{(i)}_2  \partial_i\alpha_k.
$$
Recall
the denormalization $f:\PP^1_E\ra Y$ is defined by some non-standard 
field of representatives $L=F(\omega_1+\alpha_1\epsilon,\ldots,\omega_r+\alpha_r\epsilon)$, which in turn
specifies the scalars $\alpha_i\in E$.

\begin{proposition}
\mylabel{differential equation}
In the above setting, the derivation $D\in \Der_F(E[u^{-1}])$  belongs to $\lieg$ if and only if the   partial differential equations  
$$
\lambda_1 = \Delta\lambda_0,\quad \mu_1^{(k)} = \Delta\mu_0^{(k)}, \quad\mu^{(k)}_0=\alpha_k\lambda_0,\quad \quad \Psi_k=\Delta\Phi_k 
$$
hold, with   $\Delta=\sum_{i=1}^r\alpha_i\partial_i$ from residue class ring \eqref{quotient weyl algebra} of the Weyl algebra.
\end{proposition}

\proof
The first two equations express the respective conditions that
$$
D(u^3)=\lambda_0+\lambda_1u+\ldots+\lambda_4u^4\quadand D(\omega_ku^2)= \mu_0^{(k)}+\mu_1^{(k)}u+\mu_2^{(k)}u^2,
$$
taken modulo $u^2$, belongs to the field of representatives $L\subset E[u]/(u^2)$. The third conditions means that in the Laurent polynomial
$$
D(\omega_k+\alpha_ku) = \alpha_ku^2P(u^{-1}) + Q_k(u^{-1}) + u\sum_{i=1}^rQ_i(u^{-1})\frac{\partial\alpha_k}{\partial\omega_i},
$$
the term degree $d=-2$ vanishes. This, together with the first two equations, then implies that the above has no terms
of degree $d<0$, and is thus a polynomial in $u$. The final equation means that this polynomial, taken modulo $u^2$,
belongs to $L\subset E[u]/(u^2)$.
\qed

\medskip
The first three equations ensures that  $\lambda_0\in E$ already determines $\lambda_1,\mu_0^{(k)},\mu_1^{(k)}$.
The meaning of $\Delta\Phi_k = \Psi_k$ is less clear to me. However, the   condition simplifies  dramatically 
if we choose   the field of representatives
$L=F(\omega_1+\alpha_1\epsilon,\ldots,\omega_r+\alpha_r\epsilon)$ inside the local Artin ring $\O_{Y,b}/\maxid_b^2=E[\epsilon]$ 
in a very special way:

\begin{corollary}
\mylabel{vector field p=2}
Suppose we have  $\alpha_k=\omega_k$ for $1\leq k\leq r$. Then the derivation
$$
D= (u^{-2}+1)\frac{\partial}{\partial u} + \sum_{i=1}^r\omega_i(u^{-2}+u^{-1})\frac{\partial}{\partial\omega_i}
$$
inside $\Der_F(E[u^{-1}])$ belongs to $\lieg=H^0(Y,\Theta_{Y/F})$. Furthermore, it has $D^{[2]}=0$, satisfies $D(u^3)=1+u^2$,
and vanishes on the $\omega_k+\alpha_ku$.
\end{corollary}

\proof
For the derivation at hand we get $\lambda_0=\lambda_2=1$ and $\lambda_1=\lambda_3=\lambda_4=0$ and
$\mu_0^{(k)}=\mu_1^{(k)}=\omega_k$ and $\Phi_k=\Psi_k=0$. In turn, the partial differential equations in the proposition
are obviously fulfilled. A direct computation reveals  the stated conditions on $D$.
\qed

\medskip
Suppose now $p=3$. This is     handled in an analogous  way, we only state the result. Consider derivations $D$
with coefficients $P(u^{-1})=\lambda_0u^{-3} + \lambda_1u^{-2}+\ldots + \lambda_3$ and $Q_i(u^{-1}) = \mu^{(i)}_0 u^{-3} + \mu^{(i)}_1 u^{-2} + \ldots + \mu^{(i)}_3$,
and set 
$$
\Phi_k=\alpha_k\lambda_0 + \mu_3^{(k)} + \sum_{i=1}^r \mu_2^{(i)}\partial_i\alpha_k\quadand
\Psi_k = \alpha_k\lambda_2 + \sum_{i=1}^r \mu_3^{(i)}\partial_i \alpha_k.
$$
 
\begin{proposition}
\mylabel{differential equation p=3}
In the above setting, the derivation $D\in \Der_F(E[u^{-1}])$  belongs to $\lieg$ if and only if the   partial differential equations  
$$
\lambda_1 = \Delta\lambda_0,\quad \mu_0^{(k)} = \mu_1^{(k)}=0,\quad \mu_2^{(k)}=-\alpha_k\lambda_0,\quad \Psi_k =\Delta\Phi_k
$$
hold, with   $\Delta=\sum_{i=1}^r\alpha_i\partial_i$ from residue class ring \eqref{quotient weyl algebra} of the Weyl algebra.
\end{proposition}

Specializing the coefficients that determine   $L=F(\omega_1+\alpha_1\epsilon,\ldots,\omega_r+\alpha_r\epsilon)$ 
we get:

\begin{corollary}
\mylabel{vector field p=3}
Suppose we have  $\alpha_k=\omega_k$ for $1\leq k\leq r$. Then the derivation
$$
D= (u-u^{-1})\frac{\partial}{\partial u} + \sum_{i=1}^r\omega_i(u^{-1}-1)\frac{\partial}{\partial\omega_i}
$$
inside $\Der_F(E[u^{-1}])$ belongs to $\lieg=H^0(Y,\Theta_{Y/F})$. Furthermore, it has $D^{[3]}=D$, satisfies $D(u^2)=1-u^2$,
and vanishes on the $\omega_k+\alpha_ku$.
\end{corollary}

\section{Constructions involving   nilpotents}
\mylabel{Constructions nilpotents}

As in the previous section,  $F$ denotes an imperfect ground field of characteristic $p>0$,
and $X$ is a regular curve where $E=H^0(X,\O_X)$ is a height-one extension of  degree $[E:F]=p^r$ for some $r\geq 1$.
Again we fix an $E$-rational point $a\in X$, together with a uniformizer $u\in \O_{X,a}^\wedge$,
and consider the local Artin ring $\O_{X,a}/\maxid_a^2=E[\epsilon]$.

We now choose a subfield $L\subset E[\epsilon]$ of degree $[L:F]=p^s$ for some $0\leq s<r$,
and some $L$-vector subspace $H\epsilon\subset E\epsilon $ of codimension one. We also assume
that $\Gamma(\O_X)\cap L=F$, and now 
form the subring $\Lambda=L+H\epsilon$. Note that the sum is direct, but the first summand is not
aligned to the decomposition $E[\epsilon]=E\oplus E\epsilon$.
Set
$$
A=\Spec(E[\epsilon])=\{a\}\quadand B=\Spec(L+H\epsilon)=\{b\},
$$
and consider the resulting cocartesian square
\begin{equation}
\label{square  with nilpotents}
\begin{CD}
A   @>>>    X\\
@VVV        @VVfV\\
B   @>>>    Y.
\end{CD}
\end{equation}
This defines a new integral curve $Y$ with singular locus $\Sing(Y)=\{b\}$.
Arguing as for Proposition \ref{standard models}, one sees $h^0(\O_Y)=1$ and $h^1(\O_Y)=h^1(\O_X)+1$.
Furthermore, the local rings $\O_{Y,y}$ are unibranch and Gorenstein.

We now seek to understand the complete local ring $\O_{Y,b}^\wedge$ at the $L$-valued singularity $b\in Y$.
Recall that $[E:F]=p^r$ and $[L:F]=p^s$. 
Thus $[\kappa(a):L]=p^{r-s}$, and   $\dim_L(H)$ coincides with the integer $n=p^{r-s}-1$.
We now choose a $p$-basis  and an  $L$-basis   
$$
\omega_1+\alpha_1\epsilon,\ldots,\omega_s+\alpha_s\epsilon\in L\quadand \beta_1\epsilon,\ldots,\beta_n\epsilon\in H,
$$
all belonging to  the over-ring $E[\epsilon]$. This defines   scalars $\alpha_i,\beta_j\in E$
and   $p$-independent $\omega_i\in E$. 
The ensuing cartesian square 
$$
\begin{CD}
E[\epsilon]				@<<<	E[u]\\
@AAA					@AAA\\
F[\omega_i+\alpha_i\epsilon,\beta_j\epsilon]	@<<<	R
\end{CD}
$$
defines a ring $R$ that is finitely generated over $F$.
Note that the indices in the lower left corner 
are meant to run over $1\leq i\leq s$ and $1\leq j\leq n$, a  convention that will apply throughout.

As in Section \ref{Constructions fields representatives}, 
we consider the maximal ideal $\maxid=R\cap uE[u]$. The ensuing formal completion $\widehat{R}$
sits an an analogous cartesian square, with $E[[u]]$ instead of $E[u]$ in the top right corner,
and comes with an identification   $\O_{Y,b}^\wedge=\widehat{R}$.
The following will be a key   observation:

\begin{proposition}
\mylabel{not complete intersection}
If $p^{r-s}\geq 5$  then the ring $R$ is not locally of complete intersection.
\end{proposition}

\proof
Clearly, the  monomials $\beta_ju$ belong to $\maxid$, and are linearly independent modulo $\maxid^2$,
because this   holds for their images under the canonical map  $R/\maxid^2\ra L+H\epsilon$.
We claim that these elements generate the cotangent space $\maxid/\maxid^2=\widehat{\maxid}/\widehat{\maxid}^2$, such that $\edim(R_\maxid)=n$.
To see this, consider a formal power series $\sum_{k\geq 1}\lambda_ku^k$ belonging to $\widehat{\maxid}$.
Subtracting an $L$-linear combination of the $\beta_ju$ we achieve $\lambda_1=0$.
By Lemma \ref{field-theoretic fact} below, the canonical maps $\Sym^k_L(H\epsilon)\ra \Sym^k_E(E\epsilon)$ are surjective for all $k\geq 2$.
In particular  we can write
$$
\lambda_2\epsilon^2 = \sum P_m(\omega_1+\alpha_1\epsilon,\ldots,\omega_s+\alpha_s\epsilon)\beta_{m_1}\epsilon\cdot\beta_{m_2}\epsilon
=
\sum P_m(\omega_1,\ldots,\omega_s)\beta_{m_1}\beta_{m_2}\epsilon^2,
$$
with  multi-indices $m=(m_1,m_2)$ and  $P_m\in F[T_1,\ldots,T_s]$. Note that in the outer terms, the symbol $\epsilon^2$ signifies a    monomial
in the second symmetric power rather than a ring element. Comparing coefficients gives   $\lambda_2  = \sum P_m(\omega_1,\ldots,\omega_s)\beta_{m_1}\beta_{m_2}$,
an equation in the field $E$. This in turn yields
$$
\lambda_2 u^2 = \sum P_m(\omega_1+\alpha_1u,\ldots,\omega_r+\alpha_su)\cdot (\beta_{m_1}u)(\beta_{m_2}u) + \ldots,
$$
where the omitted terms have order $\geq 3$. 
This reduces our task to $\lambda_2=0$, and likewise we achieve $\lambda_3=0$.
Now both factors in  $\sum_{k\geq 4}\lambda_ku^k=u^2\cdot\sum_{k\geq 4}\lambda_ku^{k-2}$ belong to $\widehat{\maxid}$, so the product lies in $\widehat{\maxid}^2$.
Summing up, the $\beta_ju$ form a basis form the cotangent space $\widehat{\maxid}/\widehat{\maxid}^2$.

We now choose a section for the residue class map $R\ra R/\maxid=L$. The resulting continuous homomorphism
$$
\varphi: L[[x_1,\ldots,x_n]]\lra \widehat{R},\quad x_j\longmapsto \beta_ju
$$
is surjective, by the previous paragraph. Seeking a contradiction, we assume that $R$ is locally of complete intersection.
Then the ideal of relations 
$\ideala=\Kernel(\varphi)$ is generated by $n-1$ elements. Write $\primid=(x_1,\ldots,x_n)$ for the maximal ideal.
Then $\ideala\subset\primid^2$, so according to Lemma \ref{not complete intersection} below, we have
\begin{equation}
\label{dimension module}
\dim_L(\ideala+ \primid^3)/\primid^3\leq\edim(\widehat{R})-\dim(\widehat{R})= n-1.
\end{equation}
Recall that the canonical map $\Sym^2_L(H\epsilon)\ra \Sym^2_E(E\epsilon)\simeq E$ is surjective. Its kernel $U$ 
is thus an  $L$-vector space  of dimension 
\begin{equation} 
\label{dimension kernel}
\dim_L(U) = \binom{n-1+2}{2} - (n+1) = \frac{n^2-n-2}{2}.
\end{equation}
As explained above, each 
vector $\sum\lambda_m\beta_{m_1}\epsilon\cdot\beta_{m_2}\epsilon$ from $U$ gives an equation of  in the ring $\widehat{R}\subset E[[u]]$ of the form
$$
0=\sum P_m(\omega_1+\alpha_1u,\ldots,\omega_s+\alpha_su) \cdot(\beta_{m_1}u)(\beta_{m_2} u) + \ldots,
$$
where the missing terms have order $\geq 3$.
Thus $\dim_L(\ideala+ \primid^3)/\primid^3\geq \dim_L(U)$.
Substituting \eqref{dimension module} and  \eqref{dimension kernel} we arrive at the inequality $n-1\geq (n^2-n-2)/2$, and thus $n\leq 3$.
But $n=p^{r-s}-1\geq 4$ by assumption, giving the desired contradiction.
\qed

\medskip
We are thus only interested in the situation $p^{r-s}\leq 4$. This gives three cases,
and in particular we have $p\leq 3$ and $r-s\leq 2$.

\begin{proposition}
\mylabel{complete intersection}
Suppose $p^{r-s}\leq 4$. Then   $R$ is locally of complete intersection. Disregarding the $F$-structure, the
formal completion $\widehat{R}$ is given by the following table, for certain scalars   $\lambda,\mu,\gamma,\delta\in F$:
$$
\begin{array}[b]{lll}
\toprule
p	& \widehat{R}& s	\\
\midrule
3	& L[[x,y]]/(x^3-\lambda y^3)& r-1	\\
2	& L[[x,y]]/(x^4-\mu y^2)& r-1	\\
2	& L[[x,y,z]]/(x^2-\gamma z^2, y^2-\delta z^2)& r-2	\\
\bottomrule
\end{array}
$$
\end{proposition}

\proof
First note that   $p^{r-s}\leq 4$ allows exactly for the three combinations of $p$ and $s$ occurring in the table.
Now choose a section for the residue class projection $\widehat{R}\ra L$.

Suppose first that $p=2$ and $s=r-1$. Then $H\epsilon\subset E\epsilon$ is generated by a single element $\beta_1\epsilon$
over $L$. One easily checks that the assignment $x\mapsto \beta_1u$ and $y\mapsto u^2$ gives the desired presentation
of $\widehat{R}$, with $\lambda=\beta_1^2$.
For $p=2$ and $s=r-2$ the $L$-vector space $H\epsilon$ is three-dimensional, and one checks that $x\mapsto \beta_1u$ and
$y\mapsto \beta_2u$ and $z\mapsto \beta_3u$ leads to the desired description, with $\gamma=(\beta_3/\beta_1)^2$
and $\delta=(\beta_2/\beta_1)^2$. 
Finally, for  $p=3$ and $s=r-1$ we take $x\mapsto \beta_1u$ and $y\mapsto \beta_2u$, with $\lambda=(\beta_2/\beta_1)^3$.
\qed

\medskip
From now on, we specialize to the case that 
$$
X=\PP^1_E=\Proj E[T_0,T_1] =\Spec E[u] \cup \Spec E[u^{-1}]
$$
is the projective line over $E$, with $u=T_1/T_0$ and $a=(0:1)$. Then 
$$
Y=\Spec R \cup \Spec E[u^{-1}].
$$
By construction  the $\omega_1,\ldots,\omega_{s}\in E$ are $p$-independent. 
So we can extend them to a $p$-basis $\omega_1,\ldots,\omega_r\in E$.
The  differentials $du,d\omega_i$ form a basis for the K\"ahler differentials
on the regular locus $U=Y\smallsetminus\{b\}$, and we write $\partial/\partial u, \partial/\partial\omega_i$ for the dual
basis in the tangent sheaf $\Theta_{U/F}$.
It is now possible to calculate the Lie algebra $\lieg=H^0(Y,\Theta_{Y/F})$, with the methods of Section \ref{Lie algebra}.
For the sake of brevity, we state the relevant findings:

\begin{proposition}
\mylabel{regular nilpotent}
Suppose $p^{r-s}\leq 4$ and $\pdeg(F)\geq 2r-s$. Then for suitable choices of $L+H\epsilon\subset E[\epsilon]$,
the resulting genus-one curve $Y$ admits a  twisted form $\tilde{Y}$   that is regular, and whose Frobenius base-change
$\tilde{Y}^{(p)}$ is isomorphic to some standard model $C=C^{(i)}_{r,\Lambda,F}$.
\end{proposition}

\proof
We have to consider three cases, and start with $p=3$ and $r-s=1$.
Then $n=2$ and  $\pdeg(F)\geq r+1$. So there is a  simple height-one extension $F\subset F'$ such that $E\otimes F'$ remains a field,
and endow $T=\Spec(F')$ with the structure of an $\alpha_p$-torsor.
We choose a particular $L+H\epsilon\subset E[\epsilon]$ by setting $\alpha_i=\omega_i$ for $1\leq i\leq r-1$ and
$\beta_j=\omega_r^j$ for $1\leq j\leq 2$. Now consider the derivation
$$
D=-\omega_r^{-1}(1+u)\frac{\partial}{\partial u} +(1-u^{-1})\frac{\partial}{\partial \omega_r} + \omega_r^{-1}\sum_{i=1}^{r-1} \omega_i \frac{\partial}{\partial\omega_i}. 
$$
It belongs to $\Der_F(E[u^{-1}])$, and one directly computes
$$
D(\omega_ru) = 1\quadand D(\omega_r^2u)=\omega_ru\quadand D(\omega_k+\omega_ku) = 0 \quad (1\leq k\leq r-1).
$$
This ensures that $D$ extends to a global vector field on $Y$.
Furthermore, one computes $D^{[p]}=0$.
Consider the action of the group scheme $G=\alpha_p$ on $Y$ corresponding to $D\in H^0(Y,\Theta_{Y/F})$.
One sees that the orbit of the singular locus is Cartier, and isomorphic to the trivial $G$-torsor.
The regular twisted form $\tilde{Y}$ is constructed as in Proposition \ref{regular 1,1}.

Next suppose $p=2$ and $r-s=1$, such that $n=1$. Again $\pdeg(F)\geq r+1$, and we choose $T=\Spec(F')$ as above.
The particular $L+H\epsilon\subset E[\epsilon]$ is obtained by setting $\alpha_1=\omega_1,\ldots,\alpha_{r-1}=\omega_{r-1}$
and $\beta_1=\omega_r$, and consider the derivation
$$
D=\omega_r^{-1}(u^{-1}+u)\frac{\partial}{\partial u} + 
(u^{-2}+1)\frac{\partial}{\partial\omega_r} + 
\omega^{-1}_r\sum_{i=1}^{r-1}\omega_i( u^{-1}+1)\frac{\partial}{\partial\omega_i}.
$$
This belongs to $\Der_F(E[u^{-1}])$, and one  directly computes 
$$
D(\omega_ru)=0\quadand D(\omega_ru^2)=1+u^2 \quadand D(\omega_k + \omega_ku) =0 
$$
for all $1\leq k\leq r-1$, and furthermore $D^{[p]}=0$. The argument proceeds as in the previous paragraph.

Finally, suppose $p=2$ and $r-s=2$, which is the most challenging case. Now $n=3$ and $\pdeg(F)\geq r+2$.
We choose  $L+H\epsilon\subset E[\epsilon]$ by setting 
$$
\alpha_1=\omega_1,\ldots,\alpha_{r-2}=\omega_{r-2}\quadand \beta_1=\omega_{r-1},\quad\beta_2=\omega_r,\quad\beta_3=\omega_{r-1}\omega_r.
$$
Consider the derivation
$$
D=\omega_{r-1}^{-1}(1+u)\frac{\partial}{\partial u} + \frac{\partial}{\partial\omega_{r-1}} + \omega_r\omega_{r-1}^{-1} (u^{-1}+1)\frac{\partial}{\partial \omega_r}
+ \omega_{r-1}^{-1}\sum_{i=1}^{r-2}\omega_i\frac{\partial}{\partial \omega_i}.
$$
It belongs to $\Der_F(E[u^{-1}])$, and one immediately checks
$$
D(\omega_{r-1}u)=1\quadand  D(\omega_r u)=0\quadand  D(\omega_{r-1}\omega_ru) = \omega_ru,
$$
and also  $D(\omega_k+\omega_ku)=0$ for all $1\leq k\leq r-2$. This ensures $D\in\Gamma(Y,\Theta_{Y/F})$.
Furthermore, one computes  $D^{[2]}=0$.
By symmetry, the  same formula but with  $\omega_{r-1}$ and $\omega_r$ interchanged  defines another such
$D'\in\Gamma(Y,\Theta_{Y/F})$. One computes $[D,D']=0$, and thus gets an inclusion of
the restricted Lie algebra $k^2\subset\Gamma(Y,\Theta_{Y/F})$, where on the left both bracket and
$p$-map are trivial. This corresponds to an action of the group scheme $G=\alpha_p^{\oplus 2}$ on  $Y$.
The orbit of the singular point $b\in Y$ is the effective Cartier divisor defined by the element
$\omega_{r-1}\omega_ru\in \O_{Y,b}$. 
The desired twisted form $\tilde{Y}$ that is regular arises as follows:
Now choose two simple height-one extensions $F'$
and $F''$ so that $E\otimes F'\otimes F''$ remains a field, endow their spectra with the structure
of an $\alpha_p$-torsor,   consider the resulting diagonal $G$-torsor $T=\Spec(F'\otimes F'')$,
and set  $\tilde{Y}=T\wedge^GY$.
\qed

\medskip
I would like to point out that  fixing a particular $L+H\epsilon\subset E[\epsilon]$
and finding suitable derivations $D\in\Der_F(E[u^{-1}])$ was a long and sometimes painful process of matrix computations and guesswork,
not at all reflected in the above comparatively short arguments.

\section{Verification of some technical facts}
\mylabel{Verification}

In this section we established the facts from commutative algebra and field theory used in the previous section.
We start with the former.
Recall that a noetherian ring $R$ is called \emph{locally of complete intersection} 
if for every prime $\primid$, the corresponding complete local ring
$R_\primid^\wedge$ is isomorphic to    $A/(f_1,\ldots,f_r)$, where $A$ is some   local noetherian ring
that is complete and regular, and $f_1,\ldots,f_r\in\maxid_A$ is  a regular sequence.
Note that if  $R$ itself is local, it suffices to verify this condition with the maximal ideal $\primid=\maxid_R$,
by \cite{Avramov 1975}, Corollary 1.

Suppose that $R$ is a complete local noetherian ring, and write it as $R=A/\ideala$
for some   local noetherian ring $A$ that is complete and regular, with $\edim(R)=\edim(A)$.
Fix some $n\geq 0$ with $\ideala\subset\maxid_A^n$ and form $M=(\ideala + \maxid_A^{n+1})/\maxid_A^{n+1}$.
The latter is annihilated by $\maxid_A$, and thus becomes a vector space over the residue field $k=A/\maxid_A=R/\maxid_R$.
Intuitively speaking, all relations have order at least $n$, and the vector space $M$ measures  how many of them involve terms
of order $n$.

\begin{lemma}
\mylabel{abstract not complete intersection}
If $\dim_k(M)>\edim(R)-\dim(R)$, then the ring $R$ is not locally of complete intersection.
\end{lemma}

\proof
Seeking a contradiction, we suppose that $R$ is locally of complete intersection.
By \cite{EGA IVd}, Proposition 19.3.2 the ideal $\ideala$ is generated by a regular sequence $f_1,\ldots,f_r$.
Krull's Principal Ideal Theorem gives $r=\edim(R)-\dim(R)$, and the Nakayama Lemma ensures $ \dim_k(\ideala\otimes k)\leq r$.
  The Isomorphism Theorem gives an identification
$M=\ideala/(\ideala\cap \maxid_A^{n+1})$, and hence $\dim_k(M)\leq \dim_k(\ideala\otimes_R k)$.
Combing these inequalities we get the desired contradiction.
\qed
 
\medskip
We also   used a  purely field-theoretic fact: Let $L\subset E$ be a finite field extension, $E\epsilon$ be a one-dimensional
$E$-vector space, where $\epsilon $ denotes a basis vector, and $H\epsilon$ be some $L$-linear subspace of codimension one.

\begin{lemma}
\mylabel{field-theoretic fact}
In the above setting, the canonical maps  $\Sym^k_L(H\epsilon)\ra \Sym^k_E(E\epsilon)$, $k\geq 2$ are surjective provided that $[E:L]\geq 3$.
\end{lemma}

\proof
Without restriction the one-dimensional vector space is $E$ itself, and we write $H\subset E$  for the   $L$-linear subspace
of codimension one.
The canonical map between symmetric powers becomes the multiplication map $\Sym^k_L(H)\ra E$.
It suffices to treat the case $k=2$.
Write $H\cdot H\subset E$ for the image of the multiplication map. This has codimension at most one,
because for each non-zero $a\in H$ the multiplication map $x\mapsto ax$ is injective.
Seeking a contradiction, we assume  that $H\cdot H\subset E$ is a hyperplane.
It contains, for each non-zero $a\in H$, the subvector space $aH$. The latter is also a hyperplane,
and thus $aH=H\cdot H$.

We claim that $H\subset H\cdot H$. If not,  $H\cap (H\cdot H)\subset E$ has codimension two.
Since $[E:L]\geq 3$, the intersection contains some $a\neq 0$. We just saw that we can write it in the form
$a=ab$ for some $b\in H$. Then $b=1$, and we infer that $H=1H$ is contained in $H\cdot H$, contradiction.
This actually establishes $H=H\cdot H$. 

We then even have $1\in H$: Since $\dim_L(H)=[E:L]-1\geq 2$, we find some non-zero $a\in H$,
can write it as $a=ab$ for some $b\in H$, and infer $b=1$.
It follows that the subfield $L\subset E$ is contained in the hyperplane $H\subset E$, and thus
$H\subset E$ is an $L$-subalgebra. It must be a field, because our field extensions are finite.
Set $d=[E:F]$. The Degree Formula $[E:L]=[E:H]\cdot [H:L]$  
gives $[E:H]=d/(d-1)$. Rewriting the fraction as $1+1/(d-1)$ we conclude $d=2$ or $d=0$, contradiction.
\qed

\section{Twisted ribbons and genus-zero curves}
\mylabel{Twisted ribbons}

Our analysis of genus-one curves depends on a thorough understanding of regular genus-zero curves and
their behaviour under base-changes, and we collect the relevant observations in this section.

Let $S$ be a base scheme, $Z$ be a scheme, and $\shL$ be an invertible sheaf.
As introduced by  Bayer and Eisenbud (\cite{Bayer; Eisenbud 1995}, Section 1), a \emph{ribbon for $Z$ with respect to $\shL$} is a triple $(X,i, \varphi)$,
where $X$ is a scheme, $i:Z\ra X$ is a closed embedding corresponding to  a sheaf of ideals $\shI\subset \O_X$
with $\shI^2=0$, and $\varphi:\shL\ra\shI$ is an isomorphism respecting the $\O_Z$-module structure.
A ribbon is called \emph{split} if the inclusion $i:Z\ra X$ admits a retraction $r:X\ra Z$.
Each ribbon   gives a short exact sequence
$$
0\lra \shL\lra \Omega^1_{X/S}|Z\lra \Omega^1_{Z/S}\lra 0,
$$
which in turn defines an extension class $\alpha\in\Ext^1(\Omega^1_{Z/S},\shL)$.
The ribbons from a category in the obvious way, and $(X,i,\varphi)\mapsto \alpha$ is functorial.
As explained in loc.\ cit.\, Theorem 1.2 this functor induces a bijection between the set of isomorphism class of
ribbons and the group $\Ext^1(\Omega^1_{Z/S},\shL)$.

Now suppose that $(X,i,\varphi)$ be a ribbon on $Z$ with respect to $\shL$, with class $\alpha\in \Ext^1(\Omega^1_{Z/S},\shL)$.
Let  $D\subset Z$ be any effective Cartier divisor. Set $\shL'=\shL(D)$, and let $\alpha'\in\Ext^1(\Omega^1_{Z/S},\shL')$
be the image of $\alpha$ under the canonical map $\shL\ra\shL'$.
By the universal property, the blowing-up morphism $\Bl_D(Z)\ra Z$ is an isomorphism.
Now regard $D$ as a closed subscheme on $X$. The blowing-up $X'$ comes with a canonical morphisms
$i':Z=\Bl_D(Z)\ra\Bl_D(X)= X'$. The following  is a very useful observation,
the proof being as for \cite{Bayer; Eisenbud 1995}, Theorem 1.9:

\begin{lemma}
\mylabel{blowing-up ribbons}
In the above setting, the blowing-up $\Bl_D(X)$ is a ribbon on $Z$ with respect to the invertible sheaf $\shL'$, with class $\alpha'$.
\end{lemma}

Now let $F$ be a ground field, for the moment of arbitrary characteristic $p\geq 0$.
Clearly,  a double line $2L\subset \PP^2$ is a ribbon on the line $L=\PP^1$ with respect to $\shL=\O_{\PP^1}(-1)$.
From $\Ext^1(\Omega_L,\O_L(-1))=H^1(\PP^1,\O_{\PP^1}(1))=0$ we see that each such ribbon
is split. We write it as $X_0=\PP^1\oplus\O_{\PP^1}(-1)$. Clearly, this
is a genus-zero curve, and the interpretation as double line shows that
it is locally of complete intersection. Since $\O_{2L}(-1)$ coincides with the dualizing sheaf $\omega_{X_0}$,
we see that every twisted form $X$  is also  a quadric curve in $\PP^2$.
We need the following structural result:

\begin{theorem}
\mylabel{base-change genus-zero}
Let $X$ be a twisted form of $X_0=\PP^1\oplus\O_{\PP^1}(-1)$ that is regular, and $a\in X$ a closed
point whose residue field $E=\kappa(a)$ has degree two. Then the following holds:
\begin{enumerate}
\item The field extension $F\subset E$ is simple of height one, and in particular $p=2$.
\item The only singularity on $X\otimes E$ is the rational point corresponding to $a\in X$.
\item The normalization  of  $X\otimes E$ is isomorphic to $\PP^1_{E'}$, where $F\subset E'$ has height one and degree $[E':F]=p^2$.
\end{enumerate}
Conversely, for every simple height-one extension $F\subset F'$ the base-change $X\otimes F'$ is
singular if and only if $F'\simeq\kappa(a)$ for some point $a\in X$ as above.
\end{theorem}

\proof
(i) Since the ribbon $X_0$ is everywhere singular, the local ring $\O_{X,a}$ must be geometrically singular,
hence $F\subset E$ is not separable, according to \cite{Fanelli; Schroeer 2020}, Corollary 2.6.
From $[E:F]=2$ we deduce $p=2$, and that  the extension has height one.

(iii) By construction, the point $a'\in X\otimes E$ corresponding to $a\in X$ is rational, 
so its local ring is   singular.  
According to \cite{Schroeer 2010}, Lemma 1.3 the curve $X\otimes E$ is integral. 
Let $X'=\Bl_{a'}(X\otimes E)$ be the blowing-up with reduced center $\{a'\}$.
Using that a corresponding blowing-up of the split ribbon $X_0=\PP^1\oplus\O_{\PP^1}(-1)$ is the
split ribbon $X_0'=\PP^1\oplus\O_{\PP^1}$ by Lemma \ref{blowing-up ribbons}, we see that $X'$ is a twisted form of $X'_0$,
and in particular $h^0(\O_{X'})=2$ and $h^1(\O_{X'})=0$.
In light of \cite{Schroeer 2010},  Proposition 1.4 the field $E'=H^0(X',\O_{X'})$ is a height-one extension of $F$.
By construction, it contains $E$, and we infer $[E':F]=[E':E][E:F]=p^2$.
Let $f:X'\ra X\otimes E$ be the canonical morphism. Then $D=f^{-1}(a')$ is an effective Cartier divisor
such that the canonical map $\Gamma(\O_{X'})\ra\Gamma(\O_D)$ is bijective, as one checks for the corresponding
blowing-up $X'_0\ra X_0$. It follows that $X'\simeq \PP^1_{E'}$.

(ii) Let $A$ and $B$ be the ramification and branch locus for $f:\PP^1_{E'}\ra X\otimes E$.
Using Proposition \ref{fundamental facts} we get $h^0(\O_A)-h^0(\O_B)=1$. Thus $f_*(\O_X)/\O_Y=f_*(\O_A)/\O_B$
is supported by a single point, which must be $a'$.

Finally, suppose that $F'$ is a simple height-one extension so that $X\otimes F'$ becomes singular.
The latter is integral by \cite{Schroeer 2010}, Lemma 1.3, and we write $X'$ for its normalization.
Forming the conductor square for the normalization map $X'\ra X\otimes F'$, we get an exact sequence
$$
0\lra H^0(\O_{X\otimes F'})\lra H^0(\O_{X'})\oplus H^0(\O_A)\lra H^0(\O_B)\lra 0.
$$
The field $E'=H^0(\O_{X'})$ has height one, and we write $h^0(\O_{X'})=2^r$ and $h^0(\O_A)=d2^r$ for some $r\geq 0$ and $d\geq 1$.
Then $h^0(\O_B)=d2^{r-1}$, according to Proposition \ref{fundamental facts}. From the above exact sequence
we get $1-(2^r+d2^{r-1}) + d2^r=0$, or in other words $(2-d)2^{r-1}=1$. The only solution is $d=1$ and $r=1$.
Thus $B\subset X\otimes F'$ is the inclusion of an $F'$-valued point $b'$. The image $b\in X$ is not a rational point,
hence the inclusion $\kappa(b)\subset \kappa(b')=F'$ must be an equality.
\qed

\section{Genus-one curves with singularities}
\mylabel{Genus-one curves}

Fix a ground field $F$ of characteristic $p>0$, and let $Y$ be a genus-one curve
that is not geometrically regular. The ultimate goal is to understand the situation when $Y$ is regular.
As a preparation, we   assume here that $Y$ is integral,  geometrically unibranch, locally of complete
intersection, and that some  local rings $\O_{Y,y}$ are  singular.
We furthermore assume that the Fitting ideals for $\Omega^1_{Y/F}$ are locally free, 
and that the purely inseparable field extension $E=H^0(X,\O_X)$ resulting from the normalization map $f:X\ra Y$
has height-one.
These conditions perhaps    appears   artificial, but is precisely the situation that
happens on certain base-changes from the regular case.
 Write $[E:F]=p^r$ and form the conductor square
$$
\begin{CD}
A	@>>>	X\\
@VVV		@VVfV\\
B	@>>>	Y.
\end{CD}
$$

\begin{proposition}
\mylabel{first facts}
In the above situation, the following holds:
\begin{enumerate}
\item The scheme $X$ is a regular genus-zero curve over the field $E$.
\item The normalization map $f:X\ra Y$ is a universal homeomorphism.
\item The finite $E$-scheme $A$ has degree two.
\item For each $a\in A$, we have $h^0(\O_{A,a})=2h^0(\O_{B,f(a)})$.
\item The field $F=\Gamma(\O_Y)$ is the intersection   $\Gamma(\O_X)\cap \Gamma(\O_B)$ inside $\Gamma(\O_A)$.
\end{enumerate}
\end{proposition}

\proof
We prove (i) by contradiction: Suppose $H^1(X,\O_X)\neq 0$. By the exact sequence \eqref{conductor sequence}, 
this $E$-vector space is one-dimensional as $F$-vector space.
It follows that $E=F$, and thus $h^i(\O_Y)=h^i(\O_X)$ for all $i\geq 0$. 
Now Proposition \ref{fundamental facts} gives $h^0(\O_B)=h^0(\O_A)$, hence
the inclusion $\Gamma(\O_B)\subset \Gamma(\O_A)$ is an equality. In turn, $f:X\ra Y$ is an isomorphism,
so $Y$ is regular, contradiction. This shows that $X$ is a genus-zero curve over $E$.
Condition (ii) holds because  the local rings $\O_{Y,y}$ are unibranch by assumption.

The remaining statements also follow from Proposition \ref{fundamental facts}:
(iv) and (v) are immediate consequences.
To see (iii) we use $h^0(\O_A)=2h^0(\O_B)$ to obtain $h^0(\O_X)+h^0(\O_B) = h^0(\O_A)$.
Combining these equation we see $h^0(\O_A)=2h^0(\O_X)$,
we conclude that the  $E$-algebra $\Gamma(\O_A)$ has degree two.
\qed

\medskip
From $[\Gamma(\O_A):E]=2$ and $\Gamma(\O_B)\subset\Gamma(\O_A)$ we see that there are only the   four possibilities given by the columns of the
following table:
$$
\begin{array}{lllll}
\toprule
A	& \text{non-reduced}	& \text{non-reduced	}& \text{disconnected}	& \text{integral}	\\
\midrule
B	& \text{non-reduced}	& \text{reduced}	& \text{disconnected}	& \text{integral}	\\
\bottomrule
\end{array}
$$
We shall see that the first three possibilities lead to $X=\PP^1_E$, and also the last case can be reduced to this by 
some further base-change.

\begin{proposition}
\mylabel{non-reduced in non-reduced}
Suppose $B$ and $A$ are non-reduced.  Then 
$$
p\leq 3\quadand X=\PP^1_E\quadand \Gamma(\O_A)=E[\epsilon] \quadand \Gamma(\O_B)=L+H\epsilon,
$$
for some subfield $L\subset\Gamma(\O_E)$ and some $L$-hyperplane $H\epsilon\subset E\epsilon$.
Moreover, we have $p^{r-s}\leq 4$, where $[L:F]=p^s$.
\end{proposition}

\proof
Choose a non-zero nilpotent element $\epsilon\in \Gamma(\O_B)$. This gives a non-zero map
$E[\epsilon]\ra\Gamma(\O_A)$. Since the ring of dual numbers is a principal ideal ring,
with only three ideals $(\epsilon^i)$, $0\leq i\leq 2$, we infer that
the non-zero map is injective. It actually is bijective, because both sides 
are vector spaces over $E$ of the same dimension.  
It follows that both  $A=\{a\}$ and $B=\{b\}$ are singletons. 
The statements on $\Gamma(\O_B)$ follow from Corollary \ref{subring}.
The genus-zero curve $X$ contains an $E$-valued point, so by Proposition \ref{genus-zero and genus-one}
it is isomorphic to $\PP^1_E$. Finally, the statement $p^{r-s}\leq 4$ follows from Proposition \ref{not complete intersection},
and in particular $p\leq 3$.
\qed

\begin{proposition}
\mylabel{reduced in non-reduced}
Suppose $B$ is reduced and $A$ is non-reduced. Then
$$
p\leq 3\quadand X=\PP^1_E\quadand \Gamma(\O_A)=E[\epsilon].
$$  
Moreover, the subrings $\Gamma(\O_B)$ and $\Gamma(\O_X)$ are fields of representatives inside the local Artin ring $\Gamma(\O_A)$.
\end{proposition}

\proof
Choose a non-zero nilpotent element $\epsilon\in \Gamma(\O_A)$. As in the preceding proof we infer that $\Gamma(\O_A)=E[\epsilon]$,
that $X=\PP^1_E$, and that   both $A=\{a\}$ and $B=\{b\}$ are singletons. Now  $L=\Gamma(\O_B)$ must be  a field.
From 
$2=\length_{\O_{B,b}}(\O_{A,a}) = \dim_L E[\epsilon] = 2 [E:L]$
we conclude that $L\subset E[\epsilon]$ is a field of representatives. Obviously, the same holds for the
subfield $E=\Gamma(X,\O_X)$.
\qed

\begin{proposition}
\mylabel{reduced in disconnected}
Suppose  $A$ is disconnected. Then 
$$
p=2\quadand X=\PP^1_E\quadand \Gamma(\O_A)=E\times E\quadand \Gamma(\O_B)=L'\times L'',
$$ 
where $L',L''\subset E$ are subfields with $[E:L']=[E:L'']=2$.
\end{proposition}

\proof
Since the $E$-algebra $\Gamma(\O_A)$ has degree two, it must be isomorphic to $E\times E$.
The genus-zero curve $X$ over $E$ contains an $E$-valued point, and is thus isomorphic to $\PP^1_E$.
Obviously $\Gamma(\O_B)=L'\times L''$. The inclusion $L'\subset E$ is strict.
Using that $Y$ is Gorenstein we infer  $p^{r-s}=[E:L']=2$. 
\qed

\medskip
It remains the case that  both $B$ and $A$ are integral. We shall reduce this to the preceding three cases.
The field  $\Gamma(\O_A)=\kappa(a)$  contains   both $ \Gamma(X,\O_X)=E$ and   $\Gamma(\O_B)=\kappa(b)$.
More precisely, $E$  has degree two over both subfields $\kappa(a)$ and $\kappa(b)$,
is purely inseparable over the latter, and $\Gamma(\O_Y)=F$
is the intersection of the two.  

\begin{proposition}
\mylabel{integral in integral}
In the above situation, we have $p=2$, and there is a quadratic extension $F\subset F'$ with the following properties:
\begin{enumerate}
\item The base-changes $X\otimes F'$ and $Y\otimes F'$ remain integral.
\item The field of global sections $E'=\Gamma(\O_{X'})$ for the  normalization $X'$ of   $X\otimes F'$ has height-one
over $F'$.
\item The ramification scheme $A'$ for the normalization map  $X'\ra Y\otimes F'$ is non-integral.
\end{enumerate}
\end{proposition}

\proof
The characteristic must be two, because $\kappa(b)\subset\kappa(a)$ has degree two and is purely inseparable.
The inclusion of sets $E\cup\kappa(b)\subset\kappa(a)$ must be  strict
(\cite{A 4-7}, Chapter V, \S7, No.\ 4, Lemma 1), and we choose some $\zeta\in\kappa(a)$
from the complement. This is a generator,  over both $\kappa(a)$ and $\kappa(b)$, and we have $\zeta^2\in\kappa(b)$.
If $E\subset\kappa(a)$ is purely inseparable, $\zeta^2$ is   contained in $E$, and thus in $F=E\cap \kappa(b)$,
and we conclude that  $\xi=\zeta$ has degree two over $F$.
If $E\subset\kappa(a)$ is separable, we have $\zeta^2\not\in E$, and there is an equation $\zeta^2+\zeta+\lambda=0$ for some $\lambda\in E$.
Squaring the equation gives $\zeta^4+\zeta^2+\mu=0$, where $\mu=\lambda^2$ belongs to $F$.
Now $\xi=\zeta^2$ has degree two over $F$.
In both cases we have a quadratic extension $F'=F(\xi)$ such that the minimal polynomial of $\xi$ does not split over $E$,
and thus $E(\xi)=F'\otimes E$ stays a field.

The base-change $X\otimes_FF'=X\otimes_EE(\xi)$ remains integral, according to \cite{Schroeer 2010}, Lemma 1.3,
and the same then follows for the birational curve $Y\otimes_FF'$, which establishes (i).
Assertion (ii) is a consequence of loc.\ cit., Proposition 1.4. 
It remains to verify (iii). By construction,
the generator $\xi$ belongs to $\kappa(a)=\Gamma(\O_A)$, so the base-change $A\otimes F'$ ceases to be integral.
We are done if $X\otimes F'$ remains normal, so let us assume this is not the case.
Then $F\subset F'$ is not separable, hence a height-one extension, and $X$ is a twisted form of the double line.
According to Theorem \ref{base-change genus-zero},  
the base-change $X\otimes F'=X\otimes_EE(\xi)$ acquires a singularity, and its normalization $X'$ is a projective line over some $E'$ 
containing $\xi$ and being contained in $E^{1/p}$.
The ramification scheme for $X'\ra X\otimes F'$ is the schematic fiber over the singular point in $X\otimes F'$,
and is non-reduced. It follows that the ramification locus for $X'\ra Y\otimes F'$ is non-reduced as well.
\qed

\section{Proof of the main results}
\mylabel{Proof}

This section contains the   proofs for  our main results, which where already stated in Section \ref{Formulation}.
Let $F$ be a ground field of characteristic $p\geq 0$.

\medskip
\emph{Proof of Theorem \ref{properties genus-one}}.
Suppose that  that $Y$ is a genus-one curve that is regular
but not geometric regular. Write  $r=\edim(\O_{Y,\eta}/F)$ for its geometric generic embedding dimension.
The task is to establish assertions (i)--(iii) from the enunciation of Theorem \ref{properties genus-one}.
First note that the ground field $F$ must be imperfect (\cite{EGA IVa}, Chapter 0, Theorem 22.5.8), and therefore $p>0$.
According to \cite{Schroeer 2010}, Proposition 1.5 there is a subextension $F'\subset F^{1/p}$
such that the base-change $Y\otimes F'$ remains integral, and that  its normalization $X'$ 
has $\Gamma(X',\O_{X'})=F^{1/p}$. The scheme  $X'$ is a regular genus-zero curve over $F^{1/p}$, according to 
Proposition \ref{first facts}. 
Write $B'\subset Y\otimes F'$ and  $A'\subset X'$ for the branch scheme and the 
ramification scheme for the normalization map $X'\ra Y\otimes F'$.
the latter is an  effective Cartier divisor, and its coordinate ring has degree two over the field $F^{1/p}$, again
by Proposition \ref{first facts}.  The height-one extension $F'\subset F^{1/p}$ is finite, and we set $[F^{1/p}:F']=p^r$.

Suppose for the moment that the ramification scheme $A'$ is non-integral. Then $X'=\PP^1_{F^{1/p}}$  in light of
the results from  Section \ref{Genus-one curves}. 
More precisely, we have three possible cases, and the ensuing situations are as follows:
If both  $A'$ and $B'$ are non-reduced, then   $p\leq 3$, according to Proposition \ref{non-reduced in non-reduced}.
Moreover, in the description of the coordinate ring $\Gamma(\O_{B'})=L'+H'\epsilon$ as a subring of $\Gamma(\O_{A'})=F^{1/p}[\epsilon]$
we have $p^{r-s}\leq 4$, where $[L':F']=p^s$. It follows that 
$Y\otimes F'$ is a twisted form of the standard model 
$C^{(i)}_{r,\Lambda,F'}$ with $0\leq i\leq 2$ in characteristic two, and $0\leq i\leq 1$ in characteristic three.
If $A'$ is non-reduced and $B'$ is reduced, then again $p\leq 3$, and $Y\otimes F'$ is a twisted form 
of $C^{(i)}_{r,\Lambda,F'}$ with $i=0$, now by Proposition \ref{reduced in non-reduced}.
If $A'$ is disconnected,  then Proposition \ref{reduced in disconnected} tells us that $p=2$,
and that $Y\otimes F'$ is a twisted form of some  $C^{(1,1)}_{r,F',\Lambda}$.
 
Now suppose that $A'$ and $B'$ are  integral.   Proposition \ref{integral in integral} then tells us that  $p=2$, and that 
we find some subextension $F''\subset F^{1/p^2}$ containing $F'$ such that 
$Y\otimes F''$ remains integral, and has a normalization whose ramification scheme 
becomes non-integral. In turn, the preceding paragraph applies
with $F''$ instead of $F'$.
\qed

\medskip
\emph{Proof of Theorem \ref{existence genus-one}}.
Let $r\geq 0$ be an integer, and suppose that the characteristic and $p$-degree of the ground field $F$,
together with thy symbol $i$ are as in the table of the enunciation of the theorem.
The task is to find a twisted form of $C=C^{(i)}_{r,\Lambda,F}$ that is regular, for some suitable $\Lambda$.

If  $p=2$ and $\pdeg(F)\geq r+2$ and $i=(1,1)$ then the desired twisted form exists by Proposition \ref{regular 1,1}.
Suppose now $p=2$ and $\pdeg(F)\geq r+1$ and $0\leq i\leq 2$. We then apply Proposition \ref{regular nilpotent}.
Note that the case $i=0$ is already clear by the theory of quasi-elliptic curves,
and $i=1$ is also a consequence of Proposition  \ref{regular representative}.
For $p=3$ and $\pdeg(F)\geq r+1$ and $0\leq i\leq 1$ we again apply Proposition \ref{regular nilpotent}.
\qed


\end{document}